\documentclass[11pt]{article}
\newtheorem{thm}{Theorem}
\newtheorem{pro}{Proposition}[section]
\newtheorem{lem}{Lemma}[section]
\newtheorem{assumption}{Assumption}[section]
\newtheorem{dnt}{Definition}

\newtheorem{proof}{Proof}

\catcode`@=11 \@addtoreset{equation}{section} \catcode`@=12
\def\proof{{\noindent \bf Proof} \hspace{0.3cm}}
\def\qed{\hfill \mbox{\raggedright \rule{.07in}{.1in}} \\ \vspace{0.1cm}}
\usepackage{subfigure}
\usepackage{graphicx}
\usepackage{amssymb,epsfig,amsmath,times,cite,bbm,lineno,color}
\usepackage{psfrag}
\title{\bf Transmission dynamics of an SIS model with age structure on heterogeneous networks}
\author{Shanshan Chen$^1$,  Michael Small$^{2,3}$, Yizhou Tao$^1$, Xinchu Fu$^{1,}$\thanks{Corresponding author. Tel: +86-21-66132664; Fax: +86-21-66133292; Email addresses: xcfu@shu.edu.cn (Xinchu Fu), sschen1214@163.com (Shanshan Chen)}
\\
{\it \small $^1$Department of Mathematics, Shanghai University,
Shanghai 200444, China}\\
{\it \small $^2$School of Mathematics and Statistics, University of Western Australia, 6009, Australia}\\
{\it \small $^3$Mineral Resources, CSIRO, Kensington, 6151, Australia} }

\date{}
\begin{document}

\maketitle

\begin{abstract}
\noindent Age at infection is often an important factor in epidemic dynamics. In this paper a
disease transmission model of SIS type with age dependent infection on a heterogeneous network is discussed.
The model allows the infectious rate and the recovery rate to vary and depend on the age of the infected individual at the time of infection.
We address the threshold property of the basic reproduction number and present the global dynamical
properties of the disease-free and endemic equilibria in the model.
Finally, some numerical simulations are carried out to illustrate the main results.
The combined effects of the network structure and the age dependent factor on the disease dynamics are displayed.

\vspace{0.3cm}

\end{abstract}

\section{Introduction}

\indent
Infectious diseases remain a major challenge for human society. Epidemic diseases (cholera,
tuberculosis, SARS, influenza, Ebola virus, etc.) continue to have both a major impact on human beings and an economic cost to society. Any gain we make in understanding the dynamics and control of epidemic transmission therefore has potential for significant impact --- and hence has been the focus
of scientific research and attracted much attention~\cite{1,2}.

Epidemic dynamic models provide a theoretical method for quantitative studies of infectious diseases. Since Kermack
and Mckendrick proposed two fundamental epidemic models, the SIS and SIR compartmental models, to study disease
transmission~\cite{3,4}, epidemic models for the transmission of infectious diseases have been
studied extensively. These classical compartmental models are important tools in analyzing the
spread and control of infectious diseases, but usually neglect the population structure or assume
that all the individuals have the same possibility to contact the others --- they are most effective for well-mixed homogeneous populations with a substantial penetration of infection.  However, the spreading
of infectious diseases is primarily via specific contacts between individuals, emerging diseases start with a relative small number of infectives, and the possibility to
contact others is heterogeneous. Therefore, depicting the spread of  disease processes on a
contact network can be more realistic~\cite{5,6,7,8}. Currently, the most popular transmission models on
networks are based on mean-field approximations and follow the framework initially proposed by Pastor-Satorras
and Vespignani~\cite{9,10,11}. They were the first to study SIS and SIR epidemic
models on a scale-free network and showed that the epidemic threshold is infinitesimal
in the limit of  a large number of links and nodes. Since then, a great deal of epidemiological
research work followed on scale-free (and other) networks~\cite{12,13,14,15,16,17}.

For some epidemic diseases, such as scarlet fever, poliomyelitis
 and HFMD (hand-foot-and-mouth), the process of
their transmission is related to age and some public health and preventative policies for those diseases depend on the
age structure of host population. Hence, in order to reflect the effect of demographic behavior of
individuals, researchers have begun to examine age-structured epidemic models. The
pioneering work in age-structured epidemic models was that of Hoppensteadt~\cite{18,19}, since then,
the importance of age structure in epidemic models has been recently stressed by
many authors\cite{20,21,22,23}.
Although age-dependent epidemic models have been studied extensively, all these models were
established on homogeneous networks --- in essence a convenient approximation to the homogeneous well-mixed population.
There are still few significant results concerning age-structured epidemic models on
complex (in this case, scale-free) networks.

The main purpose of this paper is to obtain threshold results for an age-structured epidemic
model on scale-free networks. A scale-free network is characterized by a power-law degree
distribution $p(k)\sim k^{-r}$~\cite{24}, where $p(k)$ is the probability that a randomly
chosen node has degree $k$, and $r$ is a characteristic exponent whose value is usually in
the range $2 <r \leq 3$. We know that for many infectious diseases, transmission can be studied
by using the SIS model with S and I representing the susceptible and infected individuals,
respectively. Based on the SIS model with age structure, our work provides new insight into
epidemic spreading dynamics.

The organization of this paper is as follows. In Section~2, we present our age-of-infection model
and give some description and assumptions. In Section~3, we analyze the existence of equilibria
and obtain the basic reproduction number. We then present some preliminaries for the analysis of
stability, which includes asymptotic smoothness of the semi-flow generated by the system and
the uniform persistence of the system.
 The main results of this paper are given in Section~4, which include the local stability and
 global stability of the disease-free and endemic equilibria. Some numerical analysis are
 performed in Section~5. Finally, in Section~6, we give conclusions and discussions.

\section{Formulation of the model}

Consider a population with connectivity modelled as a complex network N,  where each node of N is either vacant
or occupied by one individual. In an epidemic spreading process, every node has three optional
states: vacant state, susceptible state, and, infected state~\cite{25,26}. In order to consider the
heterogeneity of contacts, we divide the population into $n$~groups.
Let $S_{k}(t)$, $I_{k}(t)$, $(k=1,2,...n)$ denote the densities of susceptible and
infected nodes (individuals) with connectivity (degree) $k$ at time $t$, respectively,
and let $I_{k}(t,\tau)$ denote the density of infected individuals with respect to
the age of infection  $\tau$ at time $t$. It is obvious that
$$
I_{k}(t)=\int^\infty_0I_{k}(t,\tau){\rm d}\tau
$$
Noting that $S_{k}(t)+\int^\infty_0I_{k}(t,\tau){\rm d}\tau=N_{k}(t)$ , which describes the total density of the individuals with degree $k$ at time $t$,
then, the density of the vacant nodes with degree $k$ is $1-S_{k}-I_{k}(t)$.

In addition, as a disease spreads, a birth event occurs at a vacant node next to a non-vacant
node at rate $b$, that is to say, the empty nodes will give birth to new individuals once one of their neighbours is occupied.
Thus, the birth process depends on the number of neighboring
individuals. All individuals die at rate $\mu$, causing the occupied
node becomes vacant. Let $\beta(\tau)$, $\gamma(\tau)$ represent infectious function and
removal function with respect to age of infection $\tau$ respectively.
Therefore, the SIS epidemic model with the age-of-infection structure on a heterogeneous
network is formulated as follows:
\begin{equation}\label{eq21}\left\{\begin{array}{l}
  \frac{dS_{k}(t)}{dt}= bk[1-N_{k}(t)]\Psi_{k}-\mu S_{k}(t)-kS_{k}(t)\int^\infty_0\beta(\tau)\Theta(t,\tau){\rm d}\tau+\int^\infty_0\gamma(\tau)I_{k}(t,\tau){\rm d}\tau\\
  \frac{\partial I_{k}(t,\tau)}{\partial t}+\frac{\partial I_{k}(t,\tau)}{\partial \tau}=-(\mu+\gamma(\tau))I_{k}(t,\tau),~~~~~~~~~~~~~~~~~~~~~~~~~~0<\tau<\infty\\
 I_{k}(t,0)=kS_{k}(t)\int^\infty_0\beta(\tau)\Theta_{k}(t,\tau){\rm d}\tau,~~~~~~~~~~~~~~~~~~~~~~~~~~~~~~~~0\leq t<\infty\\
 \end{array}\right.
\end{equation}
where
$$\Theta_{k}(t,\tau)=\sum_{i=1}^n p(i|k)\frac{\varphi(i)}{i}I_{i}(t,\tau),~~\Psi_{k}(t)=\sum_{i=1}^n p(i|k)\frac{A}{i}N_{i}(t)$$
under the following initial conditions:
\begin{equation}\label{eq22}\left\{\begin{array}{l}
S_{k}(0)=S_{k0}\geq0 ~,I_{k}(0,\tau)=I_{k0}(\tau)\in L^{1}_{+}(0,\infty),(k=1,2,... n);\\
S_{k0}+\int^\infty_0I_{k0}(\tau){\rm d}\tau=N_{k0}(k=1,2,... n)\\
\end{array}\right.
\end{equation}
where~$L^{1}_{+}(0,\infty)$~is the space of functions on $[0,\infty)$ that are nonnegative and Lebesgue integrable.

The meaning of the parameters and variables of the above model are as follows:

\begin{itemize}
  \item[$\bullet$] Let~$b$ and ~$\mu$ be positive constants denoting the  birth and natural death rates of all individuals. The additional death rate induced by the infectious disease is
not considered.
  \item[$\bullet$] $\langle k\rangle$ is the average degree of the network, i.e., $\langle k\rangle=\sum_{i=1}^n ip(i)$. For a general function $\delta(k)$ , this is defined as $\langle \delta(k)\rangle=\sum_{i=1}^n \delta(i)p(i)$.  Let $p(i|k)$ be the probability that a node of degree
$k$ is connected to a node of degree $i$. In present paper, we primarily study epidemic transmissions
on uncorrelated networks, the probability is considered independent of the connectivity of the
node from which the link is emanating. Therefore, $P(i|k)=\frac{iP(i)}{\langle k\rangle}$.
  \item[$\bullet$] $\Theta_{k}(t,\tau)$~describes the probability of a link pointing to an infected
individual of age $\tau$. We note that $\varphi(k)$ is the infectivity of nodes with degree $k$,
i.e., it denotes the average number of edges from which a node with degree $k$ can transmit the
disease. Thus, $kS_{k}(t)\int^\infty_0\beta(\tau)\Theta_{k}(t,\tau){\rm d}\tau$ represents newly infected
individuals per unit time.
  \item[$\bullet$] $\Psi_{k}=\sum_{i=1}^n p(i|k)\frac{A}{i}N_{i}(t)$ is the probability of
fertility contacts between nodes with degree $k$ and its neighbours with degree $i$. The factor
$\frac{1}{i}$ accounts for the probability that one of the neighboring individual of a vacant
node with degree $i$, will activate this
vacant node at the present time step. It is assumed that, at each time step, every individual
generates the same birth contacts $A$, here $A=1$~\cite{26}. Therefore, $bk[1-N_{k}(t)]\Psi_{k}$
represents density of new born individuals per unit time.
\end{itemize}

Next we make the following assumptions on parameters, which are thought to be biologically relevant.

\begin{assumption}
Consider the system~\eqref{eq21}, we assume that,
\begin{enumerate}
  \item $b,\mu > 0$;
  \item $\beta(\tau), \gamma(\tau)\epsilon L^{1}_{+}(0,+\infty)$, with respective essential upper bounds $\overline{\beta}$ and $\overline{\gamma}$. Furthermore, there exists a constant $\alpha>0$ such that $\beta(\tau), \gamma(\tau) \geq\alpha$ for all $\tau \geq0$;
  \item $\beta(\tau), \gamma(\tau)$ are Lipschitz continuous on $R_{+}$ with Lipschitz coefficients $M_{\beta}$ and $M_{\gamma}$, respectively;
  \item For all $a\geq0$  and any k, $I_{k0}(a)>0$ . Furthermore, $\lim\limits_{ a \rightarrow \infty}I_{k0}(a)<+\infty$.
\end{enumerate}
\end{assumption}

Let us define a functional space for system~\eqref{eq21},
$$
X=X_{1}\times X_{2}\times\ldots\times X_{n},~~X_{i}=R^{1}_{+}\times L^{1}_{+}, ~i=1,2,\cdots n.
$$
Note that $X$ is a closed subset of a Banach space, and hence is a complete metric space. The norm on $X_{k}$ is taken to be $$\|X_{k}\|=|S_{k}(t)|+\int^\infty_0|I_{k}(t,\tau)|{\rm d}\tau$$
 By applying tools from~\cite{27,28} and following from Assumption 2.1, it can be verified that the solution of system~\eqref{eq21} exists and is nonnegative for any initial conditions.~Thus, for $t\geq0$ define a continuous flow $\Xi(t)$: $ X\rightarrow X$ of system~\eqref{eq21} such that
$$\Xi(t,X_{0})=\varrho(t,X_{0})=(S_{1}(t),I_{1}(t,\tau),S_{2}(t),I_{2}(t,\tau),\ldots,S_{n}(t),I_{n}(t,\tau))$$
where $\varrho(t,X_{0})$ is the solution of the model~\eqref{eq21} with initial condition $X_{0}\in X$.


~From the model~\eqref{eq21}, we get that $N_{k}$~satisfies the following differential equation,
\begin{equation}\label{eq23}
\frac{dN_{k}(t)}{dt}=bk[1-N_{k}(t)]\Psi-\mu N_{k}(t)
\end{equation}
Let $\frac{dN_{k}(t)}{dt}=0$ ,we get $N_{k}=0$, which corresponds to the equilibrium solution of
extinction, and another solution satisfies
\begin{equation}\label{eq24}
  N_{k}=\frac{bk\Psi}{\mu+bk\Psi}
\end{equation}
Putting the above equation ~\eqref{eq24} to $\Psi$, we obtain
  $$\Psi=\frac{1}{\langle k\rangle}\sum_{i} \frac{ip(i)b\Psi}{\mu+bi\Psi}$$
Noting that
$$
f(\Psi)=1-\frac{1}{\langle k\rangle}\sum_{i} \frac{bip(i)}{\mu+bi\Psi}
$$
it is clear that $f'(\Psi)>~0$, $f(1)=1-\frac{1}{\langle k\rangle}\sum_{i} \frac{bip(i)}{\mu+bi}~>~1-\frac{1}{\langle k\rangle}\sum_{i} \frac{bip(i)}{bi}~>0$.
Thus, $f(\Psi)=0$ has a unique positive solution if and only if $f(0)=1-\frac{b}{\mu}~<~0$. That is ,when $b>\mu$, the equation~\eqref{eq23} has a unique positive solution $N_{k}=N^{*}_{k}$, which satisfies
$$N^{*}_{k}=\frac{bk\Psi^{*}}{\mu+bk\Psi*}~,~\Psi^{*}=\frac{1}{\langle k\rangle}\sum_{i=1}^n p(i)N^{*}_{i}(t)$$
Therefore, from ~\eqref{eq23} and~\cite{38}, when $b~\leq~\mu$ , there is $\lim\limits_{t \rightarrow \infty}N_{k}(t)=0$, the population becomes extinct and there is no other dynamic behaviors any more. While, when $b>\mu$, $\lim\limits_{t \rightarrow \infty}N_{k}(t)=N^{*}_{k}$. Therefore, we only consider the condition of $b>\mu$ in the following sections.

Since there are the same long-playing behaviors between the original system and the limiting system. To study the stability of system~\eqref{eq21}, we consider the limiting system under which $S_{k}(t)+\int^\infty_0I_{k}(t,\tau){\rm d}\tau=N^{*}_{k}$ as follows,
\begin{equation}\label{eq25}\left\{\begin{array}{l}
  \frac{dS_{k}(t)}{dt}= bk[1-N^{*}_{k}]\Psi^{*}_{k}-\mu S_{k}(t)-kS_{k}(t)\int^\infty_0\beta(\tau)\Theta(t,\tau){\rm d}\tau+\int^\infty_0\gamma(\tau)I_{k}(t,\tau){\rm d}\tau\\
  \frac{\partial I_{k}(t,\tau)}{\partial t}+\frac{\partial I_{k}(t,\tau)}{\partial \tau}=-(\mu+\gamma(\tau))I_{k}(t,\tau)\\
 I_{k}(t,0)=kS_{k}(t)\int^\infty_0\beta(\tau)\Theta(t,\tau){\rm d}\tau\\
 \end{array}\right.
\end{equation}

Finally, we define the state space for system~\eqref{eq25} as
$$
\Gamma=\{(S_{1}(t),I_{1}(t,\cdot),\cdots ,S_{n}(t),I_{n}(t,\cdot))\in X:0~\leq~S_{k}(t)+\int^\infty_0I_{k}(t,\tau){\rm d}\tau\leq~1,k=1,\cdots, n\}
$$

 The following proposition shows that $\Gamma$ is positively invariant with respect to system~\eqref{eq25} for $\Xi$.
\begin{pro}
 $\Gamma$ is positively invariant of system~\eqref{eq25} for $\Xi$ . Moreover, $\Xi$ is point dissipative (there exists a bounded set $\Gamma\subset X$ which attracts all points in $X$). That is ,$~\forall ~t\geq0 $~, for any the solution of system~\eqref{eq25}, noting as $\Phi(t,X_{0})$ ,with the initial condition ~$X_{0}\in~\Gamma$. Then ~$\Phi(t,X_{0})\in~\Gamma$.
\end{pro}
\proof Define the arbitrary initial condition $X_{0}\in\Gamma$, which satisfies \eqref{eq22}.
Thus, $\forall t>0, N_{k0}\geq0~,\Psi(0)>0$.
From~\eqref{eq23}, we get
$$\frac{d\Psi(t)}{dt}=(b-\mu)\Psi(t)-\frac{b\Psi(t)}{\langle k\rangle}\Sigma_{i}ip(i)N_{i}(t)$$
It is obvious that $\Psi(t)=\Psi(0)e^{b-\mu-\frac{b}{\langle k\rangle}\Sigma_{i}ip(i)N_{i}(t)}~>~0$.
It follows from~\eqref{eq24} that $ N_{k}\geq 0$, and
 $$0\leq~S_{k}(t)+\int^\infty_0I_{k}(t,\tau){\rm d}\tau=N_{k}=\frac{bk\Psi}{\mu+bk\Psi}~\leq~1$$
 Therefore, $\Xi$ is point dissipative and $\Gamma$ attracts all points in X . This completes the proof. \qed

According to Assumption~2.1 and the above results, we have the following proposition.
\begin{pro}
There exists a constant $\Lambda$ satisfied $\Lambda> 1 $, then the following statements hold true for $\forall t\geq0 $ and k(k=1,2,...n):\\
(1) $0<S_{k}(t), \int^{+\infty}_{0}I_{k}(t,\tau){\rm d}\tau<\Lambda$;\\
(2) $\int^{+\infty}_{0}\beta(\tau)I_{k}(t,\tau){\rm d}\tau\leq\overline{\beta}\Lambda, \int^{+\infty}_{0}\gamma(\tau)I_{k}(t,\tau){\rm d}\tau\leq\overline{\gamma}\Lambda$, and $I_{k}(t,0)\leq \frac{\langle \varphi(k)\rangle}{\langle k \rangle}\overline{\beta}n\Lambda^{2}.$\\
(3) The function $\int^{\infty}_{0}\beta(\tau)I_{k}(t,\tau){\rm d}\tau $ and $\int^{\infty}_{0}\gamma(\tau)I_{k}(t,\tau){\rm d}\tau$ are Lipschitz continuous  with coefficient $M_{\beta}$, $M_{\gamma}$ on $R_{+}$.
\end{pro}


\section{Preliminaries}

\subsection{ Equilibria and the basic reproduction number}

Firstly, with the above boundary conditions  and initial conditions \eqref{eq22} , we obtain
$I(t,\tau)$  by integrating the second differential equation in ~\eqref{eq25} along the characteristic line
$t-\tau=constant$~\cite{29},
\begin{equation}\label{eq31}
I_{k}(t,\tau)=\left\{
\begin{array}{ll}
I_{k}(t-\tau,0)H(\tau), ~~~~t\geq\tau ; \\
I_{k}(0,\tau-t)\frac{H(\tau)}{H(\tau-t)} ~~0<t<\tau.\\
\end{array}\right.
\end{equation}
Where $H(\tau)=e^{-\mu\tau-\int^\tau_0\gamma(\xi){\rm d}\xi}$.\\
 Next, we consider the steady states of system~\eqref{eq25}. For simplicity, we note $Z_{k}(t)=I_{k}(t,0),(k=1,2,...n)$. System~\eqref{eq25} always has a disease-free equilibrium~$E_{0}: I^{0}_{k}=0, S^{0}_{k}=\frac{bk[1-N^{*}_{k}]\Psi^{*}}{\mu}~(k=1,2...n)$.

Then, let us investigate the positive equilibrium of system~\eqref{eq25}. Any positive equilibrium $E^{*}: S^{*}_{k},I^{*}_{k}(\tau)(k=1,2...n)$  should satisfy the following equations,
\begin{equation}\label{eq32}\left\{\begin{array}{l}
   bk[1-N^{*}_{k}]\Psi^{*}-\mu S^{*}_{k}-kS^{*}_{k}\int^\infty_0\beta(\tau)\Theta^{*}(\tau){\rm d}\tau+\int^\infty_0\gamma(\tau)I^{*}_{k}(\tau){\rm d}\tau=0\\
 \frac{ dI^{*}_{k}(\tau)}{d\tau}=-(\mu+\gamma)I^{*}_{k}(\tau)\\
  Z^{*}_{k}=kS^{*}_{k}\int^\infty_0\beta(\tau)\Theta^{*}(\tau){\rm d}\tau\\
  \Theta^{*}(\tau)=\frac{1}{\langle k\rangle}\sum_{i=1}^n \varphi(i)p(i)I^{*}_{i}(\tau)\\
\end{array}\right.
\end{equation}
For ease of notation, let$$K_{1}(\lambda)=\int^\infty_0\beta(\tau)e^{-\lambda\tau}H(\tau){\rm d}\tau,~  ~~K_{2}(\lambda)=\int^\infty_0\gamma(\tau)e^{-\lambda\tau}H(\tau){\rm d}\tau$$
We will get $E^{*}$, which satisfies
$$I^{*}_{k}(\tau)=I^{*}_{k}(0)H(\tau)=Z^{*}_{k}H(\tau)$$
 $$S^{*}_{k}=\frac{\langle k\rangle Z^{*}_{k}}{kK_{1}(0)\sum \varphi(i)p(i)Z^{*}_{i}},~~Z^{*}_{k}=\frac{ bk[1-N^{*}_{k}]\Psi^{*}-\mu S^{*}_{k}}{1-K_{2}(0)}$$

To make sure that $S^{*}_{k}>0$,~and $I^{*}_{k}(\tau)>0$ if and only if $Z^{*}_{k}>0$.
 It is clear that
 $$
 K_{2}(0)=\int^\infty_0\gamma(\tau)e^{-\mu\tau-\int^\tau_0\gamma(\xi){\rm d}\xi}{\rm d}\tau \leq \int^\infty_0\gamma(\tau)e^{-\int^\tau_0\gamma(\xi){\rm d}\xi}{\rm d}\tau = 1-e^{-\int^\infty_0\gamma(\xi){\rm d}\xi} < 1
 $$
Hence, we obtain
$Z^{*}_{k}>0 \Longleftrightarrow \frac{ bk[1-N^{*}_{k}]\Psi^{*}-\mu S^{*}_{k}}{1-K_{2}(0)}>0 \Longleftrightarrow  bk[1-N^{*}_{k}]\Psi^{*}-\mu S^{*}_{k}>0\Longleftrightarrow~\mu (S^{0}_{k}-S^{*}_{k})>0\Longleftrightarrow S^{0}_{k}-S^{*}_{k}>0~\Longleftrightarrow S^{0}_{k}>\frac{\langle k\rangle Z^{*}_{k}}{kK_{1}(0)\sum \varphi(i)p(i)Z^{*}_{i}}\Longleftrightarrow \frac{1}{\langle k\rangle}\sum k\varphi(k)p(k)S^{0}_{k}K_{1}(0)>1$.

From the above analysis, we get the following theorem.

\begin{thm}~~~
Define the basic reproduction number as follows,~$$R_{0}=\frac{K_{1}(0)}{<k>}\sum\limits^{n}_{i=1}i\varphi(i)p(i)S^{0}_{i}$$
If $R_{0}<1$,~the system ~\eqref{eq25} has a unique disease-free equilibrium $E^{0}$; if $R_{0}>1$,
there exist two equilibria $E^{0}$ and $E^{*}$, which satisfy
$S^{*}_{k}=\frac{\langle k\rangle Z^{*}_{k}}{kK_{1}(0)\sum \varphi(i)p(i)Z^{*}_{i}},I^{*}_{k}(\tau)=Z^{*}_{k}H(\tau),~~Z^{*}_{k}=\frac{ bk[1-N^{*}_{k}]\Psi^{*}-\mu S^{*}_{k}}{1-K_{2}(0)}, k=1,2,\ldots n$.
\end{thm}


\subsection{Asymptotic smoothness}

In order to prove the global stability of model~\eqref{eq25}, we need to make the following preparations. First, we establish asymptotic smoothness of the semigroup $\Xi(t)$. The semigroup $\Xi(t)$ is asymptotically smooth, if, for any
nonempty, closed and bounded set $B\subset X$ for which $\Xi$(t,B)$\subset$ B, there is a compact set $J\subset B$ such that $J$ attracts $B$.
In order to obtain it, we will need the following lemmas and proposition.
 \begin{lem}[\cite{30}]~
 For each~$t>0$, suppose $\Xi(t)=\Psi(t)+\Upsilon(t):\Gamma\rightarrow\Gamma$ has the property that $\Psi(t)$ is completely continuous
and there is a continuous function $k$:$R^{+}\times R^{+}\rightarrow R^{+}$:
such that $k(t,r)\rightarrow0$ as~$t\rightarrow\infty$ and $\|\Upsilon(t)X_{0}\| \leq k(t,r)$ if $\|X_{0}\|< r$.
Then $\Xi(t), t\geq0$, is asymptotically smooth.
 \end{lem}
 \begin{lem}[\cite{31}]
 Let $K\subset~L^{p}_{+}(0,+\infty)$ be closed and bounded where $p\geq1$. Then $K$ is compact
 if and only if the following conditions hold:\\
(i) $\lim_{h\rightarrow0}\int^{\infty}_{0}|u(z+h)-u(z)|^{p}{\rm d}z=0$ uniformly for $u\in K$. ($u(z+h)=0$ if $z+h<0$).\\
(ii) $\lim_{h\rightarrow\infty}\int^{\infty}_{h}|u(z)|^{p}{\rm d}z=0$ uniformly for $u\in K$.
 \end{lem}

\begin{pro}
Let $J\subset R$. For $j=1,2$, suppose that $f_{j}: J\rightarrow R$ is a bounded Lipschitz continuous function with bound $K_{j}$ and Lipschitz coefficient $M_{j}$. Then the product function $f_{1}\cdot f_{2}$ is Lipschitz continuous with coefficient $K_{1}M_{2}+K_{2}M_{1}$.
\end{pro}

From the above two lemmas, we have the following theorem.

 \begin{thm}
The semigroup $\Xi$ is asymptotically smooth.
\end{thm}
\proof  To apply Lemma 3.1, we define the projection of $\Xi(t)$ about any bounded set of $X$ by decomposing $\Xi(t)$  into the following two operators,
$$\Xi(t)=U(t)+C(t)$$
where, $$U(t)=(0,U_{1}(t),0,U_{2}(t),\ldots,0,U_{n}(t)),~~ C(t)=(S_{1}(t),C_{1}(t),S_{2}(t),C_{2}(t),\ldots,S_{n}(t),C_{n}(t))$$
\begin{equation}\label{eq33}
 U_{i}(t)=\left\{
\begin{array}{ll}
  I_{i}(0,\tau-t)\frac{H(\tau)}{H(\tau-t)} ~~~\tau\geq t\\
    0~~~~~~~~~~~~~~~~~~~~~~~~~~~~~~~\tau<t
\end{array}\right.,
 ~~C_{i}(t)=\left\{
\begin{array}{ll}
  Z_{i}(t-\tau)H(\tau) ~~~\tau<t\\
    0~~~~~~~~~~~~~~~~~~~~~~~~~~\tau\geq t
\end{array}\right.
\end{equation}
From equation~\eqref{eq31}, it is easy to get $\Xi(t)=U(t)+C(t)$. Then,
\begin{eqnarray*}
 \|U_{i}(t)\|&=&\int^{\infty}_{0}|U_{i}(t)|{\rm d}\tau=\int^{\infty}_{t} I_{i}(0,\tau-t)\frac{H(\tau)}{H(\tau-t)}{\rm d}\tau\\
&=& \int^{\infty}_{t} I_{i}(0,\tau-t)e^{-\mu t-\int^{\tau}_{\tau-t} \gamma(\xi){\rm d}\tau}{\rm d}\tau\\
&\leq& e^{-(\mu+\alpha) t}\int^{\infty}_{t} I_{i}(0,\tau-t){\rm d}\tau\\
&\leq& e^{-(\mu+\alpha) t}\|I_{i}(0,\tau)\|
\end{eqnarray*}
If $\|X_{0}\|<r$, we note $k(t,r)= re^{-(\mu+\alpha) t}$, then, $k(t,r) \rightarrow 0 $ as $t\rightarrow\infty$ and $\|U_{i}(t)\|\leq k(t,r)$ for any $i$.

Next, we verify that $C(t)$ is completely continuous. We need to pay more attention to the state space, since $L^{1}_{+}(0,+\infty)$ is a component of our state space $X$. Hence a notion of compactness in $L^{1}_{+}(0,+\infty)$ is necessary. In an infinite dimensional Banach space, boundedness does not
necessarily imply precompactness. Hence, we need to prove it by applying Lemma3.2.

 Suppose that $B\subset X$ is bounded for any initial condition $X_{0}\subset B$. From Proposition 2.1, it is easy to see that $S_{k}(t),(k=1,2,..n)$ remains in the compact set $[0,\Lambda]$. Thus, we only need to verify that the following conditions valid for $C_{i}(t)(i=1,2,..n)$ remaining in a precompact subset of  $L_{+}(0,+\infty)$ .

  To check condition $(ii)$ , according to~\eqref{eq33} ,
\begin{eqnarray*}
 C_{i}(t,\tau)&=&iS_{i}(t-\tau)\int^\infty_0\beta(\tau)\frac{1}{\langle k\rangle}\sum_{k=1}^n \varphi(k)p(k)I_{k}(t-\tau,\tau){\rm d}\tau H(\tau)\\
 &=&iS_{i}(t-\tau)\frac{1}{\langle k\rangle}\sum_{i=1}^n \varphi(k)p(k)\int^\infty_0\beta(\tau)I_{k}(t-\tau,\tau){\rm d}\tau H(\tau)\\
 &\leq& i\frac{\langle\varphi(k)\rangle}{\langle k\rangle}\overline{\beta}\Lambda^{2}e^{-(\mu+\alpha)\tau}
\end{eqnarray*}
 Note that for all $X_{0}\subset B$,  $\lim_{h\rightarrow\infty}\int^{\infty}_{h}|C_{i}(t,\tau)|_{X_{0}}^{p}{\rm d}\tau=0$. Therefore, $(ii)$ is satisfied for the set $C_{i}(t,B)\subset L^{1}_{+}(0,+\infty)$.\\
To check condition (i), for sufficiently small $h\in(0,t)$, we observe
\begin{eqnarray*}
 ~~&&\int^{\infty}_{0}|C_{i}(t,\tau)-C_{i}(t,\tau+h)|_{X_{0}}{\rm d}\tau\\
 &=& \int^{t-h}_{0}|Z_{i}(t-\tau-h)H(\tau+h)-Z_{i}(t-\tau)H(\tau)|{\rm d}\tau+ \int^{t}_{t-h}|Z_{i}(t-\tau)H(\tau)|{\rm d}\tau\\
 &\leq& \int^{t-h}_{0}Z_{i}(t-\tau-h)|H(\tau+h)- H(\tau)|{\rm d}\tau+ \int^{t-h}_{0}|Z_{i}(t-\tau-h)-Z_{i}(t-\tau)| H(\tau){\rm d}\tau\\
 &+&\int^{t}_{t-h}H(\tau)iS_{i}(t-\tau)\int^{\infty}_{0}\beta(\tau)\frac{1}{\langle k\rangle}\sum_{k}\varphi(k)p(k)I_{k}(t-\tau,\tau)|{\rm d}\tau
 \end{eqnarray*}
 It is clear that $H(\tau)=e^{-\mu\tau-\int^\tau_0\gamma(\xi){\rm d}\xi}\leq e^{-(\mu+\alpha)\tau} < 1$, and  $H(\tau)$ is a decreasing function.
 \begin{eqnarray*}
\int^{t-h}_{0}|H(\tau+h)- H(\tau)|{\rm d}\tau&=&\int^{t-h}_{0}H(\tau){\rm d}\tau- \int^{t-h}_{0}H(\tau+h){\rm d}\tau\\
&=&\int^{t-h}_{0}H(\tau){\rm d}\tau- \int^{t}_{h}H(\tau){\rm d}\tau\\
&=&\int^{h}_{0}H(\tau){\rm d}\tau- \int^{t}_{t-h}H(\tau){\rm d}\tau \leq h
 \end{eqnarray*}
Then, we note
$$\Delta=\int^{t-h}_{0} H(\tau)|Z_{i}(t-\tau-h)-Z_{i}(t-\tau)|{\rm d}\tau $$
 From Proposition 2.2 and Proposition 3.1, we have
 \begin{eqnarray*}
\Delta&=&\int^{t-h}_{0} H(\tau)|Z_{i}(t-\tau-h)-Z_{i}(t-\tau)|{\rm d}\tau\\
&=&i\sum_{k}\frac{\varphi(k)p(k)}{\langle k\rangle}\int^{t-h}_{0}H(\tau)|S_{i}(t-\tau-h)\int^{\infty}_{0}\beta(\tau)I_{k}(t-\tau-h,\tau){\rm d}\tau\\
&-&S_{i}(t-\tau)\int^{\infty}_{0}\beta(\tau)I_{k}(t-\tau,\tau){\rm d}\tau|{\rm d}\tau\\
&\leq& i\frac{\langle \varphi(k)\rangle}{\langle k\rangle}\frac{[M_{\beta}\Lambda+\overline{\beta}\Lambda M_{s}]h}{\mu+\alpha}
 \end{eqnarray*}
 $M_{s}$ is the coefficient of Lipschitz continuous function $S_{k}(t)(k=1,2,...n)$. It is easy to see that  $\frac{dS_{k}(t)}{dt}$ is bounded, and therefore, $S_{k}(t)$ is Lipschitz on $[0,\Lambda]$ with coefficient $M_{s}$.
$$
\int^{\infty}_{0}|C_{i}(t,\tau)-C_{i}(t,\tau+h)|_{X_{0}}{\rm d}\tau \leq \Lambda h+i\frac{\langle \varphi(k)\rangle}{\langle k\rangle}\frac{[M_{\beta}\Lambda+\overline{\beta}\Lambda M_{s}]h}{\mu+\alpha}+i\frac{\langle\varphi(k)\rangle}{\langle k\rangle}\overline{\beta}\Lambda^{2}h
$$
 This converges uniformly to 0 as $h\rightarrow 0$. Therefore, the condition (i) is verified for $C_{i}(t,X_{0})~(i=1,2...n)$. From Lemma~3.2, we have that $C_{i}(t,X_{0})$ is completely continuous. Finally, according to Lemma 3.1, we conclude that $\Xi(t,X_{0})$ is asymptotically smooth. This completes the proof. \qed

Next, we show that the solution semigroup $\Xi$ has a global compact attractor $A$ in $\Gamma$.

We first give the following definition of global attractors.
 \begin{dnt}[\cite{30}]
A set $A$ in $\Gamma$ is defined to be an attractor if $A$ is non-empty, compact and invariant, and there exists some open neighborhood
$W$ of $A$ in  $\Gamma$ such that $A$ attracts $W$. A global attractor is defined to be an attractor which attracts every point in $\Gamma$.
 \end{dnt}
 From the results above, we can get the existence of a global attractor  by applying  the following Lemma.
 \begin{lem}[\cite{30}]~
If $\Xi$ is asymptotically smooth and point dissipative in $\Gamma$, and  orbits of bounded sets are bounded
in $\Gamma$, then there is a global attractor $A$ in $\Gamma$.
\end{lem}
Propositions 2.1 and Theorem 3.2 show that the semigroup $\Xi(t)$ of system ~\eqref{eq25} is asymptotically smooth and point
dissipative on the state space $\Gamma$. The proof of Proposition~2.1 can verify that every forward
orbit of bounded sets is bounded in $\Gamma$. Therefore, by Lemma~3.3, we have the following theorem.
\begin{thm}
 The semigroup $\Xi$ generated by the system~\eqref{eq25} on the state space $\Gamma$ has a global attractor $A$ in $\Gamma$.
\end{thm}
\subsection{Uniform Persistence}
In this section we study the uniform persistence of system~\eqref{eq25}. Let us define a function $\rho:X \rightarrow R$ that $\rho=(\rho_{1},\rho_{2},...,\rho_{n})$ as
$$\rho_{k}(X)=Z_{k}(t),~~k=1,2,\ldots n$$
Before introducing the  result of persistence, we introduce the following important lemmas.
\begin{lem}[Fatou's Lemma]~~Let ${f_{n}}$ be a non-negative measurable function sequence, then it satisfies
\begin{equation}
\int\lim\limits_{t\rightarrow+\infty}\inf f_{n}\leq\lim\limits_{t\rightarrow+\infty}\inf\int f_{n}\leq \lim\limits_{t\rightarrow+\infty}\sup\int f_{n} \leq \int\lim\limits_{t\rightarrow+\infty}\sup f_{n}
\end{equation}
\end{lem}
\begin{lem}
[Fluctuation Lemma]~~Let ~$$\lim\limits_{t \rightarrow +\infty}\sup\Phi(t)=\Phi^{\infty},~\lim\limits_{t \rightarrow +\infty}\inf\Phi(t)=\Phi_{\infty}$$
and~$\Phi(t)$ be a bounded
and continuously differentiable function. Then there exist sequences ${f_{n}}$
and ${g_{n}}$ such that ${f_{n}}\rightarrow\infty$,~${g_{n}}\rightarrow\infty$, $\Phi(f_{n})\rightarrow \Phi^{\infty},~\Phi(g_{n})\rightarrow \Phi_{\infty}$.
$\Phi'(f_{n})\rightarrow 0 ,\Phi'(g_{n})\rightarrow~0$, as $n\rightarrow\infty$.
\end{lem}
\begin{lem}
If $R_{0}>1$, then there exists a positive constant $\varepsilon>0$, such that for any $k$,
\begin{equation}\label{eq36}
\limsup_{t\rightarrow\infty}Z_{k}(t)>\varepsilon
\end{equation}
\end{lem}
\proof If  $R_{0}>1$, there exists a sufficiently small $\varepsilon>0$ such that
\begin{equation}\label{eq37}
 \frac{1}{<k>}\sum\limits^{n}_{i=1}i\varphi(i)p(i)\frac{bi[1-N^{*}_{i}]\Psi^{*}-\varepsilon}{\mu}\int^\infty_0 \beta(\tau)H(\tau){\rm d}\tau >1
\end{equation}
We now show that this small $\varepsilon$ is the $\varepsilon$ in~\eqref{eq36}. We will do this by contradiction. Assume that there exists a constant $T>0$ which is sufficiently large such that
 $$Z_{k}(t)\leq\varepsilon ~~~~~~~~~~for~t \geq~T$$
Together with~\eqref{eq25}, we have
\begin{eqnarray*}
 \frac{dS_{k}(t)}{dt}&=&bk[1-N^{*}_{k}]\Psi^{*}-\mu S_{k}(t)-Z_{k}(t)+\int^\infty_0\gamma(\tau)I_{k}(t,\tau){\rm d}\tau\\
 &\geq&bk[1-N^{*}_{k}(t)]\Psi^{*}-\mu S_{k}-\varepsilon
\end{eqnarray*}
Then , according to the comparison principle,
\begin{equation}\label{eq37b}
S_{k}(t)\geq\frac{bk[1-N^{*}_{k}(t)]\Psi^{*}-\varepsilon}{\mu} ~~~~~~~~~for~t \geq~T
\end{equation}
Furthermore, it follows from ~\eqref{eq31} together with~\eqref{eq37b} and the Fluctuation Lemma, if ${g_{n}}$ is a sequence such that $Z_{k}(g_{n})\rightarrow Z_{k\infty}$, we have
\begin{equation}\label{eq38}
Z_{k}(g_{n})\geq k\frac{bk[1-N^{*}_{k}]\Psi^{*}-\varepsilon}{\mu}\int^{\infty}_{0}\beta(\tau)\frac{1}{\langle k\rangle}\sum_{i=1}^n \varphi(i)p(i)Z_{i}(g_{n}-\tau)H(\tau){\rm d}\tau
\end{equation}
Then, ~\eqref{eq38} becomes
$$1\geq\frac{1}{<k>}\sum\limits^{n}_{i=1}i\varphi(i)p(i)\frac{bi[1-N^{*}_{i}]\Psi^{*}-\varepsilon}{\mu}\int^\infty_0 \beta(\tau)H(\tau){\rm d}\tau$$
which contradicts to ~\eqref{eq37}. That is to say, the system ~\eqref{eq25} is  uniformly weakly $\rho$-~persistent. The proof is therefore complete. \qed

Suppose a total trajectory of $\Xi$ in space X is a function $\eta$: $R\rightarrow X$, such that $f(s+t)=\Xi(t,f(s))$, for $t\geq0, s\in R$.
Then, by applying Theorem 3.3, Lemma 3.5 and~\cite{32}, we have the following result.
\begin{thm}
If $R_{0}>1$, $\Xi(t)$ is uniformly strongly $\rho$-persistent, that is, there exists a sufficiently small $\varepsilon>0$ such that $\liminf_{t\rightarrow\infty}Z_{k}(t)>\varepsilon$.
\end{thm}

Finally, we introduce the result for uniform persistence of system~\eqref{eq25}.

\begin{thm}
If $R_{0}>1$, $\Xi(t)$ is uniform persistence, then, there exists a constant  $\varepsilon>0$ such that for any initial condition $X_{0}\in X$ and any $k$,
$$\liminf_{t\rightarrow\infty}\|I_{k}(t,\tau)\|\geq\varepsilon,~~    \liminf_{t\rightarrow\infty}S_{k}(t)\geq\varepsilon$$
\end{thm}
\proof  In fact, for any $k$, ($k=1,2,\ldots, n$), following from ~\eqref{eq22} and ~\eqref{eq31}, $I_{k}(0)\in L^{1}_{+}(0,\infty)$, we have
 \begin{eqnarray*}
 \|I_{k}(t,\tau)\|_{L^{1}}&=&\int^{t}_{0}I_{k}(t,\tau){\rm d}\tau +\int^{\infty}_{t}I_{k}(t,\tau){\rm d}\tau\\
 &\geq&\int^{t}_{0}Z_{k}(t-\tau)H(\tau){\rm d}\tau
\end{eqnarray*}
According to Lemma 3.4 and Theorem 3.4, we obtain that there exists a sufficiently small $\varepsilon>0$ satisfying the following inequality:
 \begin{eqnarray*}
 \liminf_{t\rightarrow\infty}\|I_{k}(t,\tau)\|_{L^{1}}&\geq&\int^{\infty}_{0}\liminf_{t\rightarrow\infty}Z_{k}(t-\tau)H(\tau){\rm d}\tau\\
 &>&\varepsilon\int^{\infty}_{0}H(\tau){\rm d}\tau=\frac{\varepsilon}{\mu+\overline{\gamma}}\triangleq\varepsilon_{1}
\end{eqnarray*}
Then, by a similar argument, we have
 \begin{eqnarray*}
\frac{dS_{k}(t)}{dt}&=&bk[1-N^{*}_{k}]\Psi^{*}_{k}-\mu S_{k}(t)-kS_{k}(t)\int^\infty_0\beta(\tau)\Theta(t,\tau){\rm d}\tau+\int^\infty_0\gamma(\tau)I_{k}(t,\tau){\rm d}\tau\\
&\geq&bk[1-N^{*}_{k}]\Psi^{*}_{k}-(\mu+k\int^\infty_0\beta(\tau)\Theta(t,\tau){\rm d}\tau )S_{k}(t)+\alpha(N^{*}_{k}-S_{k}(t))
\end{eqnarray*}
Hence, by the comparison principle, we obtain $ \liminf_{t\rightarrow\infty}S_{k}(t)\geq\frac{bk[1-N^{*}_{k}]\Psi^{*}_{k}+\alpha N^{*}_{k}}{(\mu+\alpha)+k\frac{\langle \varphi(k)\rangle}{\langle k\rangle} \overline{\beta}\Lambda}\triangleq\varepsilon_{2}$.
Therefore, we take $\varepsilon=\max\{\varepsilon_{1},\varepsilon_{2}\}$, the proof is completed.  \qed

\section{The  epidemic threshold}

\subsection{Local stability}

In this section, we first evaluate the local stability of equilibria of system~\eqref{eq25}. For convenience, we apply~\eqref{eq31} to get the following system of equations for model~\eqref{eq25}:
\begin{equation}\left\{\begin{array}{l}\label{41}
  \frac{dS_{k}(t)}{dt}= bk[1-N^{*}_{k}]\Psi^{*}-\mu S_{k}(t)-kS_{k}(t)\int^t_0\beta(\tau)\Theta_{1}(t,\tau){\rm d}\tau-B^{1}_{k}(t)+\int^t_0\gamma(\tau)Z_{k}(t-\tau)H(\tau){\rm d}\tau+B^{2}_{k}(t)\\
  Z_{k}(t)=kS_{k}(t)\int^t_0\beta(\tau)\Theta_{1}(t,\tau){\rm d}\tau+B^{1}_{k}(t)\\
  \Theta_{1}(t,\tau)=\frac{1}{\langle k\rangle}\sum_{i=1}^n \varphi(i)p(i)Z_{i}(t-\tau)H(\tau),~~t\geq\tau \\
  \Theta_{2}(t,\tau)=\frac{1}{\langle k\rangle}\sum_{i=1}^n \varphi(i)p(i)I_{i}(0,\tau-t)\frac{H(\tau)}{H(\tau-t)},~~t<\tau
\end{array}\right.
\end{equation}
with equation ~\eqref{eq31}, where,
 $$B^{1}_{k}(t)=kS_{k}(t)\int^\infty_t\beta(\tau)\Theta_{2}(t,\tau){\rm d}\tau, ~~ B^{2}_{k}(t)=\int^\infty_t\gamma(\tau)I_{k}(0,\tau-t)\frac{H(\tau)}{H(\tau-t)}{\rm d}\tau$$
It is obvious that $\lim\limits_{t \rightarrow \infty} B^{2}_{k}(t)=0, ~\lim\limits_{t \rightarrow \infty} B^{1}_{k}(t)=0$.

Then considering the limiting system associated with ~\eqref{eq31}:
\begin{equation}\left\{\begin{array}{l}\label{eq42}
  \frac{dS_{k}(t)}{dt}= bk[1-N^{*}_{k}]\Psi^{*}-\mu S_{k}(t)-kS_{k}(t)\int^\infty_0\beta(\tau)\Theta_{1}(t,\tau){\rm d}\tau+\int^\infty_0\gamma(\tau)Z_{k}(t-\tau)H(\tau){\rm d}\tau\\
  Z_{k}(t)=kS_{k}(t)\int^\infty_0\beta(\tau)\Theta_{1}(t,\tau){\rm d}\tau
\end{array}\right.
\end{equation}
 with the same initial conditions with $(2.1)$ and $I_{k}(t,\tau)$ can get from ~\eqref{eq31}, where
 $$\Theta_{1}(t,\tau)=\frac{1}{\langle k\rangle}\sum_{i=1}^n \varphi(i)p(i)Z_{i}(t-\tau)H(\tau)$$
By using the Jacobian matrix and its characteristic equation, we have the following theorem:
\begin{thm}~~If $R_{0}<1$,~ the disease-free
equilibrium $E_0$  is locally asymptotically
stable; ~and it is unstable while~$R_{0}>1$.
\end{thm}
\proof First,  linearizing \eqref{eq42} near $E_{0}$ by denoting the perturbation variables $S_{k}(t)=\tilde{S}_{k}(t)+S^{0}_{k}, Z_{k}(t)=\tilde{Z}_{k}(t), \Theta_{1}(t)=\theta(t)$,
 we obtain the following system
\begin{equation}\left\{\begin{array}{l}\label{eq43}
  \frac{d\tilde{S}_{k}(t)}{dt}= -\mu\tilde{S}_{k}(t)-kS^{0}_{k}\int^\infty_0\beta(\tau)\theta(t,\tau){\rm d}\tau+\int^\infty_0\gamma(\tau)\tilde{Z}_{k}(t-\tau)H(\tau){\rm d}\tau\\
  \tilde{Z}_{k}(t)=k\tilde{S}_{k}(t)\int^\infty_0\beta(\tau)\theta(t,\tau){\rm d}\tau\\
  \theta(t,\tau)=\frac{1}{\langle k\rangle}\sum_{i=1}^n \varphi(i)p(i)\tilde{Z}_{i}(t-\tau)H(\tau)
\end{array}\right.
\end{equation}
Analyzing the local asymptotic stability near $E_{0}$, let $\tilde{S}_{k}(t)=s_{k0}e^{\lambda t},\tilde{Z}_{k}(t)=z_{k0}e^{\lambda t}$ and substitute them into~\eqref{eq43}, we get
\begin{equation}\left\{\begin{array}{l}\label{eq44}
(\lambda+\mu)s_{k0}+\frac{kS^{0}_{k}}{<k>}K_{1}(\lambda)\sum\limits^{n}_{i=1}\varphi(i)p(i)z_{i0}-K_{2}(\lambda)z_{k0}=0\\
z_{k0}-\frac{kS^{0}_{k}}{<k>}K_{1}(\lambda)\sum\limits^{M}_{i=1}\varphi(i)p(i)z_{i0}=0\\
\end{array}\right.
\end{equation}
Then, we can get the characteristic equation of ~\eqref{eq42},
$$
\left|
\begin{smallmatrix}
\lambda+\mu&0&\cdots&0&\frac{S^{0}_{1}K_{1}(\lambda)}{<k>}\varphi(1)p(1)-K_{2}(\lambda)&\frac{S^{0}_{1}K_{1}(\lambda)}{<k>}\varphi(2)p(2)&\cdots&\frac{S^{0}_{1}K_{1}(\lambda)}{<k>}\varphi(n)p(n)\\
  0&\lambda+\mu&\cdots&0&\frac{2S^{0}_{2}K_{1}(\lambda)}{<k>}\varphi(1)p(1)&\frac{2S^{0}_{2}K_{1}(\lambda)}{<k>}\varphi(2)p(2)-K_{2}(\lambda)&\cdots&\frac{2S^{0}_{2}K_{1}(\lambda)}{<k>}\varphi(n)p(n)\\
  \vdots&\vdots&\cdots&\vdots&\vdots&\vdots&\cdots&\vdots\\
  0&0&\cdots&\lambda+\mu&\frac{nS^{0}_{n}K_{1}(\lambda)}{<k>}\varphi(1)p(1)&\frac{nS^{0}_{n}K_{1}(\lambda)}{<k>}\varphi(2)p(2)&\cdots&\frac{nS^{0}_{n}K_{1}(\lambda)}{<k>}\varphi(n)p(n)-K_{2}(\lambda)\\
  0&0&\cdots&0&1-\frac{S^{0}_{1}K_{1}(\lambda)}{<k>}\varphi(1)p(1)&-\frac{S^{0}_{1}K_{1}(\lambda)}{<k>}\varphi(2)p(2)&\cdots&-\frac{S^{0}_{1}K_{1}(\lambda)}{<k>}\varphi(n)p(n)\\
  0&0&\cdots&0&-\frac{2S^{0}_{2}K_{1}(\lambda)}{<k>}\varphi(1)p(1)&1-\frac{2S^{0}_{2}K_{1}(\lambda)}{<k>}\varphi(2)p(2)&\cdots&-\frac{2S^{0}_{2}K_{1}(\lambda)}{<k>}\varphi(n)p(n)\\
  \vdots&\vdots&\cdots&\vdots&\vdots&\vdots&\cdots&\vdots\\
  0&0&\cdots&0&-\frac{nS^{0}_{n}K_{1}(\lambda)}{<k>}\varphi(1)p(1)&-\frac{nS^{0}_{n}K_{1}(\lambda)}{<k>}\varphi(2)p(2)&\cdots&1-\frac{nS^{0}_{n}K_{1}(\lambda)}{<k>}\varphi(n)p(n)\\
  \end{smallmatrix}
       \right|=0
$$
which is equivalent to the following form,
$$\left|
\begin{matrix}
(\lambda+\mu)E_{n}&A\\
0&B\\
  \end{matrix}
       \right|=0
$$
That is to say, $(\lambda+\mu)E_{n}B=0$, therefore, the eigenvalues are $\lambda_{k}=-\mu (k=1,2,\cdots,n),$ and satisfy $B=0$,
where,
$$
\left|
\begin{matrix}
1-\frac{S^{0}_{1}K_{1}(\lambda)}{<k>}\varphi(1)p(1)&-\frac{S^{0}_{1}K_{1}(\lambda)}{<k>}\varphi(2)p(2)&\cdots&-\frac{S^{0}_{1}K_{1}(\lambda)}{<k>}\varphi(n)p(n)\\
 -\frac{2S^{0}_{2}K_{1}(\lambda)}{<k>}\varphi(1)p(1)&1-\frac{2S^{0}_{2}K_{1}(\lambda)}{<k>}\varphi(2)p(2)&\cdots&-\frac{2S^{0}_{2}K_{1}(\lambda)}{<k>}\varphi(n)p(n)\\
  \vdots&\vdots&\cdots&\vdots\\
  -\frac{nS^{0}_{n}K_{1}(\lambda)}{<k>}\varphi(1)p(1)&-\frac{nS^{0}_{n}K_{1}(\lambda)}{<k>}\varphi(2)p(2)&\cdots&1-\frac{nS^{0}_{n}K_{1}(\lambda)}{<k>}\varphi(n)p(n)\\
  \end{matrix}
       \right|=0
$$
That is,
\begin{equation}\label{eq45}
\frac{K_{1}(\lambda)}{<k>}\sum\limits^{n}_{i=1}i\varphi(i)p(i)S^{0}_{i}=1
\end{equation}
We denote ~$G(\lambda)=\frac{K_{1}(\lambda)}{<k>}\sum\limits^{n}_{i=1}i\varphi(i)p(i)S^{0}_{i}$, which satisfies $G(\lambda)$=1.\\
Suppose that $\lambda=a+bi$ is the solution of $G(\lambda)=1$ and $a\geq0$. When $R_{0}<1$,
\begin{eqnarray*}
|G(\lambda)|&=&|\frac{1}{<k>}\sum\limits^{n}_{i=1}S^{0}_{i}i\varphi(i)p(i)\int^\infty_0\beta(\tau)H(\tau)e^{-(a+bi)\tau}{\rm d}\tau|\\
&=& |\frac{1}{<k>}\sum\limits^{n}_{i=1}S^{0}_{i}i\varphi(i)p(i)\int^\infty_0\beta(\tau)H(\tau)e^{-a\tau}(\cos(b\tau)-i\sin(b\tau)){\rm d}\tau|\\
&\leq& \frac{1}{<k>}\sum\limits^{M}_{i=1}S^{0}_{i}i\varphi(i)p(i)\int^\infty_0\beta(\tau)H(\tau)e^{-a\tau}{\rm d}\tau\\
&\leq&\frac{1}{<k>}\sum\limits^{M}_{i=1}S^{0}_{i}i\varphi(i)p(i)\int^\infty_0\beta(\tau)H(\tau){\rm d}\tau=R_{0}<1.
\end{eqnarray*}
It is obvious that $|G(\lambda)|<1$ is contradictory to $G(\lambda)=1$.~Thus,~$G(\lambda)=1$ doesn't have positive solutions. So all roots of the characteristic equation are negative. Therefore, this means that if $R_{0}<1$, the disease-free equilibrium $E_{0}$ is locally asymptotically stable.

On the other hand,~$K_{1}(\lambda)=\int^\infty_0\beta(\tau)H(\tau)e^{-\lambda\tau}{\rm d}\tau$, $\lim\limits_{\lambda \rightarrow-\infty}K_{1}(\lambda)=\infty$, $\lim\limits_{\lambda \rightarrow +\infty}K_{1}(\lambda)=0$. Therefore, when $R_{0}>1$,~$G(\lambda)=1$ has positive real part, that is to say,~$E_{0}$~is unstable. This completes the proof. \qed

Next, we discuss the local stability of the positive equilibrium. We have the following Theorem.
\begin{thm}~~~If $R_{0}>1$, system \eqref{eq42} has a unique positive equilibrium point $E^{*}=(S_{k}^{*},Z_{k}^{*}),~k=1,2,...,n$, and it is locally stable.
\end{thm}
\proof  The existence of the positive  equilibrium has been obtained in Section~3. As the same method applied in the discussion of the stability of  disease-free equilibrium. First, linearizing the system \eqref{eq42} near $E^{*}$, that is, $S_{k}(t)=S^{*}_{k}+\tilde{S}_{k}(t)$, $Z_{k}(t)=Z^{*}_{k}+\tilde{Z}_{k}(t)$, we get the following linear system
\begin{equation}\left\{\begin{array}{l}\label{eq46}
  \frac{d\tilde{S}_{k}(t)}{dt}= -\mu \tilde{S}_{k}(t)-kS^{*}_{k}\int^\infty_0\beta(\tau)\theta(t,\tau){\rm d}\tau-k\tilde{S}_{k}(t)\int^\infty_0\beta(\tau)\theta^{*}(t,\tau){\rm d}\tau +\int^\infty_0\gamma(\tau)\tilde{Z}_{k}(t-\tau)H(\tau){\rm d}\tau\\
  \tilde{Z}_{k}(t)= kS^{*}_{k}\int^\infty_0\beta(\tau)\theta(t,\tau){\rm d}+k\tilde{S}_{k}(t)\int^\infty_0\beta(\tau)\theta^{*}(t,\tau){\rm d}\tau\\
  \theta^{*}(t,\tau) =\frac{1}{\langle k\rangle}\sum_{i=1}^n \varphi(i)p(i)Z^{*}_{i}(t-\tau)H(\tau)
\end{array}\right.
\end{equation}
Let $\tilde{S}_{k}(t)=s_{0k}e^{\lambda t}, \tilde{Z}_{k}(t)=z_{0k}e^{\lambda t}$, where $ s_{0k},z_{0k}$ can be determined later, and substitute them into the system~\eqref{eq46}, we obtain the following equation:
\begin{equation}\left\{\begin{array}{l}\label{eq47}
(\lambda+\mu+kbc)s_{0k}+kS^{*}_{k}C(\lambda)\sum\limits^{n}_{i=1}\varphi(i)p(i)z_{0i}-K_{2}(\lambda)z_{0k}=0\\
kbcs_{0k}+kS^{*}_{k}C(\lambda)\sum\limits^{n}_{i=1}\varphi(i)p(i)z_{0i}-z_{0k}=0
\end{array}\right.
\end{equation}
where $b=\sum\limits^{n}_{i=1}\varphi(i)p(i)Z^{*}_{i},~c=\frac{K_{1}(0)}{\langle k\rangle},~C(\lambda)=\frac{K_{1}(\lambda)}{\langle k\rangle},~\varphi(i)p(i)=\varphi_{i}P_{i}$.
Then, we analyze the local stability of the positive equilibrium $E^{*}$. From \eqref{eq47}, we can get the characteristic equation as follows:
$$\left|
\begin{smallmatrix}
\lambda+\mu+bc&0&\ldots&0&S^{*}_{1}C(\lambda)\varphi_{1}P_{1}-K_{2}(\lambda)&S^{*}_{1}C(\lambda)\varphi_{2}P_{2}&\ldots&S^{*}_{1}C(\lambda)\varphi_{n}P_{n}\\
0&\lambda+\mu+2bc&\ldots&0&2S^{*}_{2}C(\lambda)\varphi_{1}P_{1}&2S^{*}_{2}C(\lambda)\varphi_{2}P_{2}-K_{2}(\lambda)&\ldots&2S^{*}_{2}C(\lambda)\varphi_{n}P_{n}\\
\vdots&\vdots&\ldots&\vdots&\vdots&\vdots&\ldots&\vdots\\
0&0&\ldots&\lambda+\mu+nbc&nS^{*}_{n}C(\lambda)\varphi_{1}P_{1}&nS^{*}_{n}C(\lambda)\varphi_{2}P_{2}&\ldots&nS^{*}_{n}C(\lambda)\varphi_{n}P_{n}-K_{2}(\lambda)\\
bc&0&\ldots&0&S^{*}_{1}C(\lambda)\varphi_{1}P_{1}-1&S^{*}_{1}C(\lambda)\varphi_{2}P_{2}&\ldots&S^{*}_{1}C(\lambda)\varphi_{n}P_{n}\\
0&2bc&\ldots&0&2S^{*}_{2}C(\lambda)\varphi_{1}P_{1}&2S^{*}_{2}C(\lambda)\varphi_{2}P_{2}-1&\ldots&2S^{*}_{2}C(\lambda)\varphi_{n}P_{n}\\
\vdots&\vdots&\ldots&\vdots&\vdots&\vdots&\ldots&\vdots\\
0&0&\ldots&nbc&nS^{*}_{n}C(\lambda)\varphi_{1}P_{1}&nS^{*}_{n}C(\lambda)\varphi_{2}P_{2}&\ldots&nS^{*}_{n}C(\lambda)\varphi_{n}P_{n}-1\\
 \end{smallmatrix}
       \right|=0
$$
Note that $m_{i}=\frac{ibc}{\lambda+\mu+ibc},~i=1,2...n$,
we obtain
\begin{equation}\label{eq48}
\sum\limits^{n}_{i=1}[\frac{(1-m_{i})S^{*}_{i}i\varphi_{i}P_{i}C(\lambda)}{m_{i}K_{2}(\lambda)-1}+1]\cdot\textstyle\prod\limits^{n}_{i=1}(m_{i}K_{2}(\lambda)-1)=0
\end{equation}

\textbf{Case 1.}  If $\textstyle\prod\limits^{n}_{i=1}(m_{i}K_{2}(\lambda)-1)=0$, that is to say,  $K_{2}(\lambda)=1+\frac{\lambda+\mu}{ibc},i=1,2,...n$.
We assume that ~~$\lambda\geq0$, then it is obvious that
$$K_{2}(\lambda)=\int^\infty_0\gamma(\tau)e^{-\lambda\tau}H(\tau){\rm d}\tau=\int^\infty_0\gamma(\tau)e^{-(\lambda+\mu)\tau-\int^\tau_0 \gamma(\xi){\rm d}\tau}{\rm d}\tau~\leq~\int^\infty_0\gamma(\tau)e^{-\int^\infty_0 \gamma(\xi){\rm d}\tau}{\rm d}\tau~<~1$$
 However, $1+\frac{\lambda+\mu}{ibc}\geq1$, thus, the assumption is contradictory. Therefore, if $\textstyle\prod\limits^{n}_{i=1}(m_{i}K_{2}(\lambda)-1)=0$, the eigenvalue is negative.

\textbf{Case 2.}  If $\sum\limits^{n}_{i=1}[\frac{(1-m_{i})S^{*}_{i}i\varphi_{i}P_{i}C(\lambda)}{m_{i}K_{2}(\lambda)-1}+1]=0$,
substituting $m_{n}$,$b,c$ and $C(\lambda)$ into the above equation, we get
$$\sum\limits^{n}_{i=1}\frac{(\lambda+\mu)i\varphi(i)P(i)S^{*}_{i}}{i\sum \varphi(i)P(i)Z^{*}_{i}[1-K_{2}(\lambda)]+(\lambda+\mu)\frac{\langle K\rangle}{K_{1}(0)}}=\frac{K_{1}(0)}{K_{1}(\lambda)}$$
Due to $S^{*}_{k}=\frac{\langle k\rangle Z^{*}_{k}}{kK_{1}(0)\sum \varphi(i)p(i)Z^{*}_{i}}$, then,
$$\sum\limits^{n}_{i=1}\frac{(\lambda+\mu)\langle k\rangle \varphi(i)P(i)Z^{*}_{i}}{i[\sum  \varphi(i)P(i)Z^{*}_{i}]^{2}K_{1}(0)[1-K_{2}(\lambda)]+(\lambda+\mu)\langle k\rangle \sum_{i} \varphi(i)P(i)Z^{*}_{i}}=\frac{K_{1}(0)}{K_{1}(\lambda)}$$
It is obvious that
$$\frac{K_{1}(0)}{K_{1}(\lambda)}~<~\frac{(\lambda+\mu)\langle k\rangle \sum_{i} \varphi(i)P(i)Z^{*}_{i}}{[\sum_{i} \varphi(i)P(i)Z^{*}_{i}]^{2}K_{1}(0)[1-K_{2}(\lambda)]+(\lambda+\mu)\langle k\rangle \sum_{i} \varphi(i)P(i)Z^{*}_{i}}$$\\
Supposing $Re\lambda\geq0$,
then we have $|\frac{K_{1}(0)}{K_{1}(\lambda)}|\geq1$, whereas,
$$
\frac{(\lambda+\mu)\langle k\rangle \sum_{i} \varphi(i)P(i)Z^{*}_{i}}{[\sum_{i} \varphi(i)P(i)Z^{*}_{i}]^{2}K_{1}(0)[1-K_{2}(\lambda)]+(\lambda+\mu)\langle k\rangle \sum_{i} \varphi(i)P(i)Z^{*}_{i}} < 1
$$
Therefore, the assumption is also contradictory. So all the eigenvalues are negative.

To sum up, \eqref{eq48} do not have roots with positive real parts, therefore, all the roots of \eqref{eq48} have negative real parts. Therefore, the endemic equilibrium $E^{*}$ is locally asymptotically stable if $R_{0}>1$.  The proof is complete. \qed


\subsection{Global stability of equilibria}

 In this section, we first study the global stability of the disease-free equilibrium~$E_{0}$~by using the Fluctuation Lemma and Fatou's Lemma~\cite{33}. We have the following theorem.
\begin{thm}~~~If~$R_{0}<1$~, the disease-free equilibrium~$E_0$~of the system \eqref{eq25} is globally asymptotically stable.
\end{thm}
\proof  Because $Z_{k}(t)>0$, let $V_{k}(t)=\int^\infty_0\gamma(\tau)Z_{k}(t-\tau)H(\tau){\rm d}\tau$, which is nonnegative.
$$
\frac{dS_{k}(t)}{dt}=  bk[1-N^{*}_{k}]\Psi^{*}-\mu S_{k}(t)-Z_{k}(t)+V_{k}(t)~\leq~ bk[1-N^{*}_{k}(t)]\Psi^{*}-\mu S_{k}+V_{k}(t)
$$
According to the comparison principle and Proposition~2.2, it is easy to get the following equation:
$$
S_{k}(t)\leq S_{k0}e^{-\mu t}+\frac{ bk[1-N^{*}_{k}]\Psi^{*}}{\mu}(1-e^{-\mu t})+\int^t_0V_{k}(\xi)e^{-\mu(t-\xi)}{\rm d}\xi
$$
That is to say,  $\lim\limits_{t\rightarrow+\infty}\sup S_{k}(t)~=S^{\infty}_{k}~\leq~\frac{ bk[1-N^{*}_{k}]\Psi^{*}}{\mu}=S^{0}_{k}$.

Then, we verify $\lim\limits_{t\rightarrow+\infty} Z_{k}(t)=0~(k=1,2,...,n)$. From the system \eqref{eq42}, we know $\{Z_{k}(t)\}$ $(k=1,2,\cdots,n)$ is a measurable sequence of non-negative uniformly bounded functions. Based on Fatou's Lemma, we have
\begin{eqnarray*}
\lim\limits_{t\rightarrow+\infty}\sup Z_{k}(t)&=&\lim\limits_{t\rightarrow+\infty}\sup\frac{kS_{k}(t)}{\langle k\rangle}\int^\infty_0\beta(\tau)\sum^{n}_{i=1} \varphi(i)p(i)Z_{i}(t-\tau)H(\tau){\rm d}\tau\\
\lim\limits_{t\rightarrow+\infty}\sup Z_{k}(t)&\leq&\frac{kS^{\infty}_{k}}{\langle k\rangle}\int^\infty_0\beta(\tau)\lim\limits_{t\rightarrow+\infty}\sup\sum^{n}_{i=1} \varphi(i)p(i)Z_{i}(t-\tau)H(\tau){\rm d}\tau
\end{eqnarray*}
Then, we let the inequality be multiplied by $\varphi(k)p(k)$~and sum over $k$, we have
\begin{eqnarray*}
\varphi(k)p(k)\lim\limits_{t\rightarrow+\infty}\sup Z_{k}(t)&\leq&\frac{k\varphi(k)S^{0}_{k}p(k)}{\langle k\rangle}\int^\infty_0\beta(\tau)\sum_{i} \varphi(i)p(i)\lim\limits_{t\rightarrow+\infty}\sup Z_{i}(t)H(\tau){\rm d}\tau\\
\sum_{k} \varphi(k)p(k)\lim\limits_{t\rightarrow+\infty}\sup Z_{k}(t)&\leq&\frac{\sum_{k} k\varphi(k)S^{0}_{k}p(k)}{\langle k\rangle}\int^\infty_0\beta(\tau)\sum_{i} \varphi(i)p(i)\lim\limits_{t\rightarrow+\infty}\sup Z_{i}(t)H(\tau){\rm d}\tau\\
\sum_{k} \varphi(k)p(k)\lim\limits_{t\rightarrow+\infty}\sup Z_{k}(t)&\leq&\frac{K_{1}(0)}{<k>}\sum\limits^{n}_{k=1}S^{0}_{k}k\varphi(k)p(k)\sum_{i} \varphi(i)p(i)\lim\limits_{t\rightarrow+\infty}\sup Z_{i}(t)\\
\sum_{k} \varphi(k)p(k)\lim\limits_{t\rightarrow+\infty}\sup Z_{k}(t)&\leq&R_{0}\sum_{i} \varphi(i)p(i)\lim\limits_{t\rightarrow+\infty}\sup Z_{i}(t)
\end{eqnarray*}
Because $R_{0}<1$, so  if this inequality holds only when $\sum_{i}\varphi(i)p(i)\lim\limits_{t\rightarrow+\infty}\sup Z_{i}(t)=0$, that is,
$$
\lim\limits_{t\rightarrow+\infty}\sup Z_{k}(t)=0
$$
 On the other hand, $Z_{k}(t)_{\infty}~\geq~0$ because the  positive definiteness. Then, we get $\lim\limits_{t\rightarrow+\infty}Z_{k}(t)=0$.\\
 Moreover, Lemma 3.5 implies that there exists a sequence ${g_{n}}$, such that ${g_{n}}\rightarrow\infty$, then $\Phi(g_{n})\rightarrow \Phi_{\infty}$
and $\Phi'(g_{n})\rightarrow 0$, as $n\rightarrow\infty$. Thus,
 $$\frac{dS_{k}(g_{n})}{dt}=  bk[1-N^{*}_{k}]\Psi^{*}-\mu S_{k}(g_{n})-Z_{k}(g_{n})+\int^\infty_0\gamma(\tau)Z_{k}(g_{n}-\tau)H(\tau){\rm d}\tau$$
 Let~$n\rightarrow\infty$, then,
 $$ bk[1-N^{*}_{k}]\Psi^{*}-\mu S_{k\infty}-Z^{\infty}_{k}+\int^\infty_0\gamma(\tau)Z_{k\infty}H(\tau){\rm d}\tau~\leq~0$$
 Because of  $Z_{k\infty} ,~Z^{\infty}_{k}\rightarrow~0$, then, $\frac{ bk[1-N^{*}_{k}]\Psi^{*}}{\mu}~\leq~S_{k\infty}$.
 Thus $\frac{ bk[1-N^{*}_{k}]\Psi^{*}}{\mu}~\leq~S_{k\infty}~\leq~S^{\infty}_{k}~\leq~\frac{ bk[1-N^{*}_{k}]\Psi^{*}}{\mu}$. That is, $\lim\limits_{t\rightarrow+\infty}S_{k}(t)=\frac{ bk[1-N^{*}_{k}]\Psi^{*}}{\mu}=S^{0}_{k}~(k=1,2,...n)$.

To sum up, $(S_{k}(t),Z_{k}(t))\rightarrow~E_{0}$ in $\Gamma$ for any $k=1,2,\ldots, n$, as $t\rightarrow\infty$. The proof is therefore completed.  \qed


In the following, we verify the global stability of $E^{*}$ of the system~\eqref{eq42} by Lyapunov-LaSalle asymptotic stability theorem for the semiflow $\Xi(t)$. Now, to simplify the model~\eqref{eq25}, we let $\gamma(\tau)=\gamma$ be a constant.
Since $S_{k}(t)+\int^\infty_0I_{k}(t,\tau){\rm d}\tau=N^{*}_{k}$, hence
we have $\int^\infty_0I_{k}(t,\tau)d\tau=N^{*}_{k}-S_{k}(t)$. For convenience, we denote ~$bk[1-N^{*}_{k}(t)]\Psi^{*}+\gamma N^{*}_{k}=\Lambda_{k}$.
Then model~\eqref{eq25} can be represented as
\begin{equation}\label{eq49}\left\{\begin{array}{l}
  \frac{dS_{k}(t)}{dt}=\Lambda_{k}-(\mu+\gamma)S_{k}(t)-kS_{k}(t)\int^\infty_0\beta(\tau)\Theta(t,\tau){\rm d}\tau\\
  \frac{\partial I_{k}(t,\tau)}{\partial t}+\frac{\partial I_{k}(t,\tau)}{\partial \tau}=-(\mu+\gamma(\tau))I_{k}(t,\tau)\\
  I_{k}(t,0)=kS_{k}(t)\int^\infty_0\beta(\tau)\Theta(t,\tau){\rm d}\tau\\
  \Theta(t,\tau)=\frac{1}{\langle k\rangle}\sum_{n=1}^M \varphi(n)p(n)I_{n}(t,\tau)
\end{array}\right.
\end{equation}

\begin{thm}~~~If~$R_{0}>1$~, the endemic equilibrium~$E^{*}$~of the system \eqref{eq49} is globally asymptotically stable.
\end{thm}
\proof
 Firstly, we introduce the following important function, which is obtained from the linear combination of Volterra-type functions of the form
$$
g (x) = x-1-\ln x
$$
Obviously, $ g(x)\geq0$ for $x>0$ and $g'(x)=1-1/x$. Then, $g(x)$ has a global minimum at $x=1$ and $g(1)=0$. That is to say, $g(x)$ is nonnegative.

Next, constructing the positively definite Lyapunov functional
$V_{k}(t)=V_{S_{k}}(t)+V_{I_{k}}(t)$, where,
$$
V_{S_{k}}(t)=S^{*}_{k}(\frac{S_{k}(t)}{S^{*}_{k}}-1-\ln\frac{S_{k}(t)}{S^{*}_{k}})
$$
$$
V_{I_{k}}(t)=\frac{kS^{*}_{k}}{\langle k\rangle}\frac{\sum_{i}\varphi(i)p(i)I^{*}_{i}(0)}{I^{*}_{k}(0)}\int^\infty_0 \pi(\tau) I^{*}_{k}(\tau)(\frac{I_{k}(t,\tau)}{I^{*}_{k}(\tau)}-1-ln\frac{I_{k}(t,\tau)}{I^{*}_{k}(\tau)}){\rm d}\tau
$$
We denote $\pi(\tau)=\int^\infty_\tau \beta(\xi)e^{-(\mu+\gamma)(\xi-\tau)}{\rm d}\xi$, which satisfies
$$\pi(0)=\int^\infty_0 \beta(\xi)H(\xi){\rm d}\xi=~K_{1}(0)$$
$$\frac{d\pi(\tau)}{d\tau}=(\mu+\gamma)\pi(\tau)-\beta(\tau)$$
Because the positive equilibrium $E^{*}$ satisfies~\eqref{eq32},
we obtain $$I^{*}_{k}(\tau)=I^{*}_{k}(0)H(\tau)$$
$$\Lambda_{k}=(\mu+\gamma) S^{*}_{k}-kS^{*}_{k}\int^\infty_0\beta(\tau)\Theta^{*}(\tau){\rm d}\tau$$
$$ I^{*}_{k}(0)=\frac{kS^{*}_{k}}{\langle k\rangle}\sum_{i=1}^n \varphi(i)p(i)\int^\infty_0\beta(\tau)I^{*}_{i}(\tau){\rm d}\tau=\frac{kS^{*}_{k}}{\langle k\rangle}K_{1}(0)\sum_{i=1}^n \varphi(i)p(i)I^{*}_{i}(0)$$
To prove that the Lyapunov functional $V_{k}(t)$ is well-defined, it suffices to show that
$$
\triangle_{1}\triangleq \int^\infty_0 I^{*}_{k}(\tau)ln\frac{I_{k}(t,\tau)}{I^{*}_{k}(\tau)}{\rm d}\tau, ~~ \triangle_{2}\triangleq S^{*}_{k}ln\frac{S_{k}(t)}{S^{*}_{k}}
$$
are finite for all $t\geq0$ and all $k=1,2,\cdots,n$. It is clearly true from  Theorem~3.5, $\triangle_{2}$ is finite for all $t\geq0$ and all $k=1,2,\cdots,n$. Meanwhile, Assumption~2.1 and Theorem~3.5 ensure that $\triangle_{1}$ is finite for all $t\geq0$ and all $k=1,2,\cdots,n$.
Then, the derivative of $V_{k}$ along the solutions of \eqref{eq49} is
$$\frac{dV_{k}(t)}{dt}=\frac{dV_{S_{k}}(t)}{dt}+\frac{dV_{I_{k}}(t)}{dt}$$
where,
\begin{eqnarray*}
\frac{dV_{S_{k}}(t)}{dt}&=& (1-\frac{S^{*}_{k}}{S_{k}(t)})\frac{dS_{k}(t)}{dt}\\
&=&-\frac{(\mu+\gamma)}{S_{k}(t)}(S_{k}-S^{*}_{k})^{2}+kS^{*}_{k}\int^\infty_0\beta(\tau)\Theta^{*}(\tau){\rm d}\tau-kS_{k}\int^\infty_0\beta(\tau)\Theta(t,\tau){\rm d}\tau\\
& &-\frac{S^{*}_{k}}{S_{k}(t)}kS^{*}_{k}\int^\infty_0\beta(\tau)\Theta^{*}(\tau){\rm d}\tau+kS^{*}_{k}\int^\infty_0\beta(\tau)\Theta(t,\tau){\rm d}\tau
\end{eqnarray*}
\begin{eqnarray*}
\frac{dV_{I_{k}}(t)}{dt}&=&\frac{kS^{*}_{k}}{\langle k\rangle}\frac{\sum_{i}\varphi(i)p(i)I^{*}_{i}(0)}{I^{*}_{k}(0)}\int^\infty_0 \pi(\tau)I^{*}_{k}(\tau)(\frac{1}{I^{*}_{k}(\tau)}-\frac{1}{I_{k}(t,\tau)})\frac{\partial I_{k}(t,\tau)}{\partial t}{\rm d}\tau\\
&=&-\frac{kS^{*}_{k}}{\langle k\rangle}\frac{\sum_{i}\varphi(i)p(i)I^{*}_{i}(0)}{I^{*}_{k}(0)}\int^\infty_0 \pi(\tau)I^{*}_{k}(\tau)(\frac{1}{I^{*}_{k}(\tau)}-\frac{1}{I_{k}(t,\tau)})[\frac{\partial I_{k}(t,\tau)}{\partial \tau}+(\mu+\gamma)I_{k}(t,\tau)]{\rm d}\tau\\
&=&-\frac{kS^{*}_{k}}{\langle k\rangle}\frac{\sum_{i}\varphi(i)p(i)I^{*}_{i}(0)}{I^{*}_{k}(0)}\int^\infty_0 \pi(\tau)I^{*}_{k}(\tau) \frac{\partial}{\partial \tau} g[\frac{I_{k}(t,\tau)}{I^{*}_{k}(\tau)}]{\rm d}\tau\\
&=&-\pi(\infty)I^{*}_{k}(\infty)g[\frac{I_{k}(t,\infty)}{I^{*}_{k}(\infty)}]\frac{kS^{*}_{k}}{\langle k\rangle}\frac{\sum_{i}\varphi(i)p(i)I^{*}_{i}(0)}{I^{*}_{k}(0)}+\pi(0)\frac{kS^{*}_{k}}{\langle k\rangle}\frac{\sum_{i}\varphi(i)p(i)I^{*}_{i}(0)}{I^{*}_{k}(0)}I^{*}_{k}(0)g(\frac{I_{k}(t,0)}{I^{*}_{k}(0)})\\
& &+\frac{kS^{*}_{k}}{\langle k\rangle}\frac{\sum_{i}\varphi(i)p(i)I^{*}_{i}(0)}{I^{*}_{k}(0)}\int^\infty_0 [\frac{d\pi(\tau)}{d\tau}I^{*}_{k}(\tau)+\frac{dI^{*}_{k}(\tau)}{d\tau}\pi(\tau)]g(\frac{I_{k}(t,\tau)}{I^{*}_{k}(\tau)}){\rm d}\tau\\
&=&-\pi(\infty)I^{*}_{k}(\infty)g[\frac{I_{k}(t,\infty)}{I^{*}_{k}(\infty)}]\frac{kS^{*}_{k}}{\langle k\rangle}\frac{\sum_{i}\varphi(i)p(i)I^{*}_{i}(0)}{I^{*}_{k}(0)}+\pi(0)\frac{kS^{*}_{k}}{\langle k\rangle}\frac{\sum_{i}\varphi(i)p(i)I^{*}_{i}(0)}{I^{*}_{k}(0)}I^{*}_{k}(0)g(\frac{I_{k}(t,0)}{I^{*}_{k}(0)})\\
& &-\frac{kS^{*}_{k}}{\langle k\rangle}\frac{\sum_{i}\varphi(i)p(i)I^{*}_{i}(0)}{I^{*}_{k}(0)}\int^\infty_0\beta(\tau)I^{*}_{k}(0)H(\tau)g(\frac{I_{k}(t,\tau)}{I^{*}_{k}(\tau)}){\rm d}\tau\\
&=&-B+I_{k}(t,0)-I^{*}_{k}(0)+I^{*}_{k}(0)ln\frac{I_{k}(t,0)}{I^{*}_{k}(0)}-\frac{1}{\langle k\rangle}\sum_{i}\varphi(i)p(i)kS^{*}_{k}\int^\infty_0\beta(\tau)I^{*}_{i}(\tau)g(\frac{I_{k}(t,\tau)}{I^{*}_{k}(\tau)}){\rm d}\tau
\end{eqnarray*}
where $B=-\pi(\infty)I^{*}_{k}(\infty)g(\frac{I_{k}(t,\infty)}{I^{*}_{k}(\infty)})\frac{kS^{*}_{k}}{\langle k\rangle}\frac{\sum_{i}\varphi(i)p(i)I^{*}_{i}(0)}{I^{*}_{k}(0)}.$
Then, from the above two parts, we get
\begin{eqnarray*}
\frac{dV_{k}(t)}{dt}=&-&\frac{(\mu+\gamma)}{S_{k}(t)}(S_{k}-S^{*}_{k})^{2}-B+\int^\infty_0\beta(\tau)\frac{1}{\langle k\rangle}\sum_{i}\varphi(i)p(i)kS^{*}_{k}I^{*}_{i}(\tau)[1-\frac{S^{*}_{k}}{S_{k}(t)}+\frac{I_{i}(t,\tau)}{I^{*}_{i}(\tau)}\\
&-&\frac{S_{k}(t)I_{i}(t,\tau)I^{*}_{k}(0)}{S^{*}_{k}I^{*}_{i}(\tau)I_{k}(t,0)}-ln\frac{I_{k}(t,0)}{I^{*}_{k}(0)}-\frac{I_{k}(t,\tau)}{I^{*}_{k}(\tau)}+1+ln\frac{I_{k}(t,\tau)}{I^{*}_{k}(\tau)}]{\rm d}\tau\\
=~&-&\frac{(\mu+\gamma)}{S_{k}(t)}(S_{k}-S^{*}_{k})^{2}-B+\int^\infty_0\beta(\tau)\frac{1}{\langle k\rangle}\sum_{i}\varphi(i)p(i)kS^{*}_{k}I^{*}_{i}(\tau)[1-\frac{S^{*}_{k}}{S_{k}(t)}+ln\frac{S^{*}_{k}}{S_{k}(t)}\\
&+&1-\frac{S_{k}(t)I_{i}(t,\tau)I^{*}_{k}(0)}{S^{*}_{k}I^{*}_{i}(\tau)I_{k}(t,0)}+ln\frac{S_{k}(t)I_{i}(t,\tau)I^{*}_{k}(0)}{S^{*}_{k}I^{*}_{i}(\tau)I_{k}(t,0)}\\
&+&\frac{I_{i}(t,\tau)}{I^{*}_{i}(\tau)}-ln\frac{I_{i}(t,\tau)}{I^{*}_{i}(\tau)}-\frac{I_{k}(t,\tau)}{I^{*}_{k}(\tau)}+ln\frac{I_{k}(t,\tau)}{I^{*}_{k}(\tau)}]{\rm d}\tau\\
=~&-&\frac{(\mu+\gamma)}{S_{k}(t)}(S_{k}-S^{*}_{k})^{2}-B+\frac{1}{\langle k\rangle}\sum_{i}\varphi(i)p(i)kS^{*}_{k}I^{*}_{i}(0)\int^\infty_0\beta(\tau)H(\tau)[1-\frac{S^{*}_{k}}{S_{k}(t)}+ln\frac{S^{*}_{k}}{S_{k}(t)}+1\\
&-&\frac{S_{k}(t)I_{i}(t,\tau)I^{*}_{k}(0)}{S^{*}_{k}I^{*}_{i}(\tau)I_{k}(t,0)}+ln\frac{S_{k}(t)I_{i}(t,\tau)I^{*}_{k}(0)}{S^{*}_{k}I^{*}_{i}(\tau)I_{k}(t,0)}
+\frac{I_{i}(t,\tau)}{I^{*}_{i}(\tau)}-ln\frac{I_{i}(t,\tau)}{I^{*}_{i}(\tau)}-\frac{I_{k}(t,\tau)}{I^{*}_{k}(\tau)}+ln\frac{I_{k}(t,\tau)}{I^{*}_{k}(\tau)}]{\rm d}\tau
\end{eqnarray*}
Let $$a_{ki}(\tau)=\frac{1}{\langle k\rangle}\sum_{i}\varphi(i)p(i)kS^{*}_{k}I^{*}_{i}(0)$$
$$f(x)=1-x+\ln x$$
$$G_{k}(I_{k})=\int^\infty_0\beta(\tau)H(\tau)[-\frac{I_{k}(t,\tau)}{I^{*}_{k}(\tau)}+ln\frac{I_{k}(t,\tau)}{I^{*}_{k}(\tau)}]{\rm d}\tau$$
\begin{eqnarray*}
\Upsilon_{ki}(t,S_{k},I_{k}(t,\cdot))&=&\int^\infty_0\beta(\tau)H(\tau)[1-\frac{S^{*}_{k}}{S_{k}(t)}+ln\frac{S^{*}_{k}}{S_{k}(t)}+1-\frac{S_{k}(t)I_{i}(t,\tau)I^{*}_{k}(0)}{S^{*}_{k}I^{*}_{i}(\tau)I_{k}(t,0)}
+ln\frac{S_{k}(t)I_{i}(t,\tau)I^{*}_{k}(0)}{S^{*}_{k}I^{*}_{i}(\tau)I_{k}(t,0)}\\ &&-\frac{I_{k}(t,\tau)}{I^{*}_{k}(\tau)}+ln\frac{I_{k}(t,\tau)}{I^{*}_{k}(\tau)}-(-\frac{I_{i}(t,\tau)}{I^{*}_{i}(\tau)}+ln\frac{I_{i}(t,\tau)}{I^{*}_{i}(\tau)})]{\rm d}\tau
\end{eqnarray*}
Therefore, $$\frac{dV_{k}}{dt}\leq\sum_{i}a_{ki}\Upsilon_{ki}(t,S_{k},I_{k}(t,\cdot))$$
Furthermore,
\begin{eqnarray*}
\Upsilon_{ki}(t,S_{k},I_{k}(t,\cdot))&=&\int^\infty_0\beta(\tau)H(\tau)[f(\frac{S^{*}_{k}}{S_{k}(t)})+f(\frac{S_{k}(t)I_{i}(t,\tau)I^{*}_{k}(0)}{S^{*}_{k}I^{*}_{i}(\tau)I_{k}(t,0)}){\rm ]d}\tau+G_{k}(I_{k})-G_{i}(I_{i})\\
&\leq&G_{k}(I_{k})-G_{i}(I_{i})
\end{eqnarray*}
Finally, according to Corollary~3.3 and Theorem~3.1 in~\cite{34}, $V_{k},\Upsilon_{ki},G_{ki},a_{ki}$ satisfy the assumptions. Therefore, the function $V=\sum_{k}c_{k}V_{k}$ which is defined in Corollary~3.3 is a Lyapunov function for system~\eqref{eq49}. It is obvious that $V'(t)\leq0$ for the model~\eqref{eq49} and $(S_{1},I_{1}(t,\tau),\cdots,S_{n},I_{n}(t,\tau))\in\Gamma$. In addition, the largest invariant set for $V'(t)=0$ is $E^{*}$. Then, the positive solution of system~\eqref{eq49} is globally asymptotically stable according to  Theorem~5.3.1 and Corollary~5.3.1 in~\cite{35}. The proof is therefore completed.  \qed

\section{Simulations}

We know that the stability of the disease-free and endemic equilibria depends on the basic reproduction number. Here we present numerical simulations to explore the effects of various parameters on the basic reproduction number, and to support the analytic results obtained in the previous sections.

To simulate the process of system~\eqref{eq21}, we first adopt a first-order upwind scheme with forward Euler time step to process the second PDEs of~\eqref{eq21}; Consequently, we use a Runge-Kutta scheme to deal with the first ODEs of~\eqref{eq21}. The partition of mesh is $\Delta t=0.1$ and $\Delta\tau=0.2$, which satisfies the Courant-Friedrichs-Lewy stability restriction condition for PDEs.

Our simulations are based on scale-free uncorrelated networks with the degree distribution $p(k)=ck^{-r}$, where $r=2.4$, and the constant $c$ satisfies $\Sigma^{n}_{k=1}k^{-r}=1$. We set the maximum degree $n=40$. The initial values $S_{k}(0)=0.6$, $I_{k}(0,\tau)=\frac{1}{\sqrt{2\pi}}e^{\frac{-(\tau+1.15)^{2}}{2}}$ for any degree $k$.  The other parameters are chosen as $b=0.07$ and $\mu=0.06$.

First, Fig.~1 depicts the influence of $\varphi(k)$ on $R_{0}$ where $\varphi(k)$ is the infectivity of the disease. This function has many potential forms, such as $\varphi(k)=k$~\cite{9,11}, which means that the number of the  contacts  per unit time is equal to the node's degree $k$; $\varphi(k)=h$~\cite{36}, which means that the number of the  contacts  per unit time is equal to the constant $h$ and has no relation with the node's degree; or (for example), $\varphi(k)=\frac{\omega k^{a}}{1+\nu k^{a}}$~\cite{37}. Here, we discuss out results assuming this general form.
Fig.~1 shows that $R_{0}$ is monotonically increasing as $\alpha$ increases.
The bigger  $a$ is and the smaller $b$ is, the bigger $R_{0}$ is. This implies that the contact method can influence disease transition --- as expected.

Next, we present the evolution of $I_{k}(t)$ over time under different parameters from Fig.~2 and Fig.~3 with  the infectivity $\varphi(k)=k$, where $I_{k}(t)=\int^\infty_0 I_{k}(t,\tau){\rm d}\tau$.
Fig.~2 shows that the evolutions of the density of infection with regard to time and degree under the parameters chosen as $\beta(\tau)=\frac{\tau(200-\tau)}{15000},\gamma(\tau)=\frac{1}{1+\tau}$, which ensures that $R_{0}<1$.
Fig.~2(a) shows that the time series of $I_{10},~I_{20}$, $I_{30}$ and $I_{40}$. Fig.~2(b) depicts the overall trend of infection with all degrees. They clearly demonstrate that when $R_{0}<1$  the disease will gradually die out, and the disease-free equilibrium is globally asymptotically stable.  That is, $\lim\limits_{t\rightarrow+\infty} I_{k}(t)=0$. In addition, we conclude that the larger the degree, the higher the peak of infection.

In Fig.~3, we show the evolution of $I_{k}(t)$ over time with infectivity $\varphi(k)=k$ under the parameters chosen as $\beta(\tau)=\frac{\tau(200-\tau)}{15000}, \gamma(\tau)=\frac{1}{1+10\tau}$ , the basic reproduction number satisfies $R_{0}>1$. In this case, simulation indicates that the disease eventually become endemic, and tends to the endemic equilibrium which is globally asymptotically stable. Moreover, the larger the degree, the higher the endemic level.
Therefore, Fig.~2 and Fig.~3 can support the results obtained in previous sections, and ensure the global stability when $\gamma(\tau)$ is not an constant.

\begin{figure}[htbp]
\centering
\includegraphics[height=5.0cm,width=7.0cm]{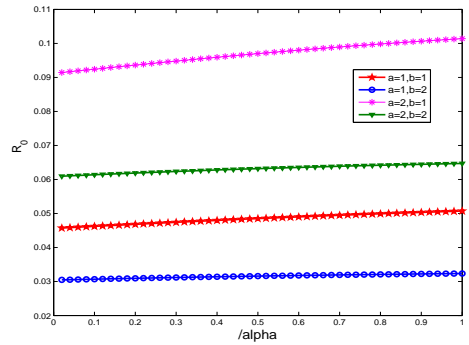}
\caption{\footnotesize{(Color online) The influence of $\varphi(k)$ on $R_{0}$, here $\varphi(k)=\frac{ak^{\alpha}}{1+bk^{\alpha}}$}}
\label{Fig1}
\end{figure}

\begin{figure}[htbp]
\centering
\subfigure[Dynamics of $I_{k}$(t) subject to time $t$]{\includegraphics[width=6.3cm]{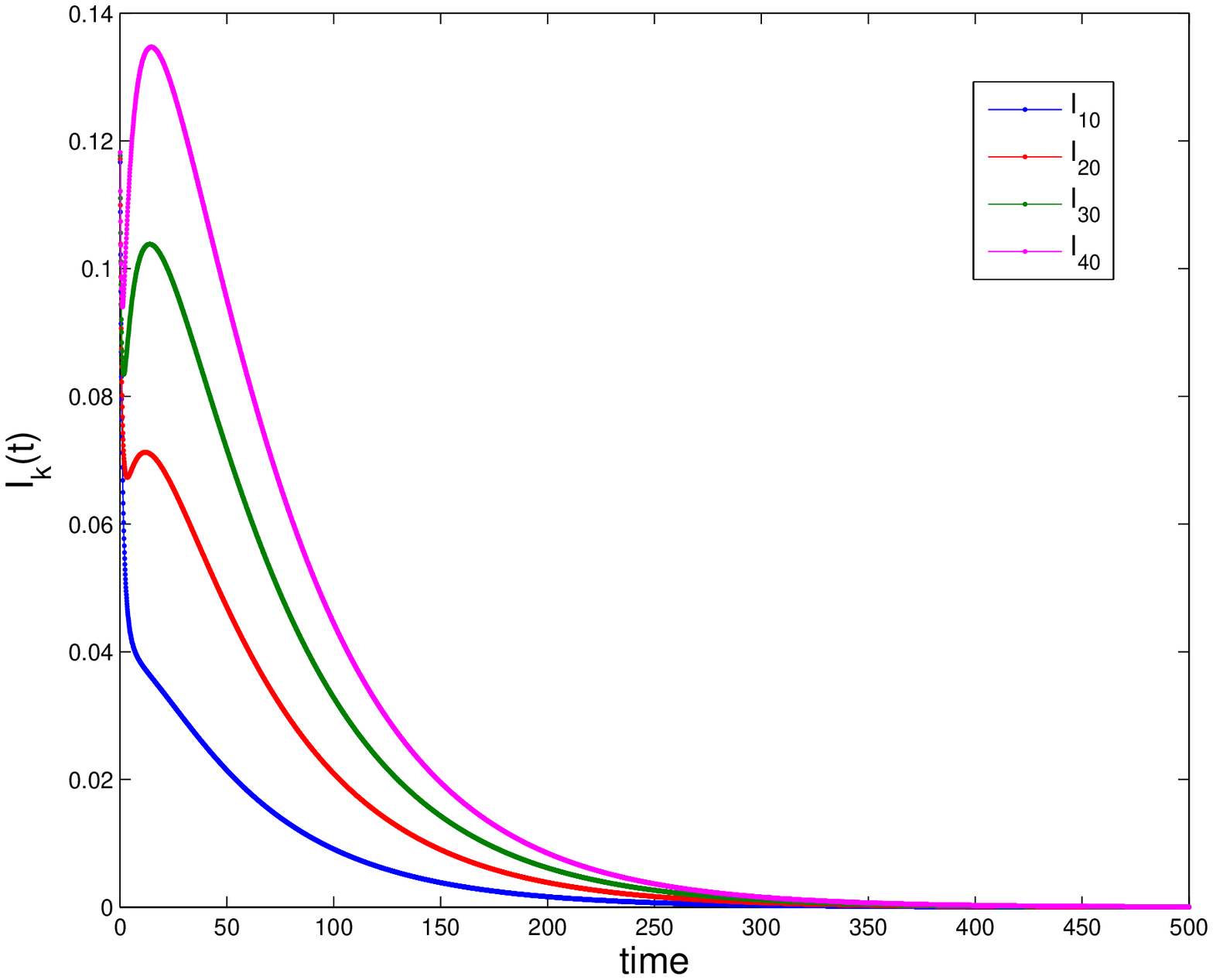}
\label{Fig2(a)}}
\subfigure[Dynamics of $I_{k}$(t) subject to time $t$]{\includegraphics[width=7.2cm]{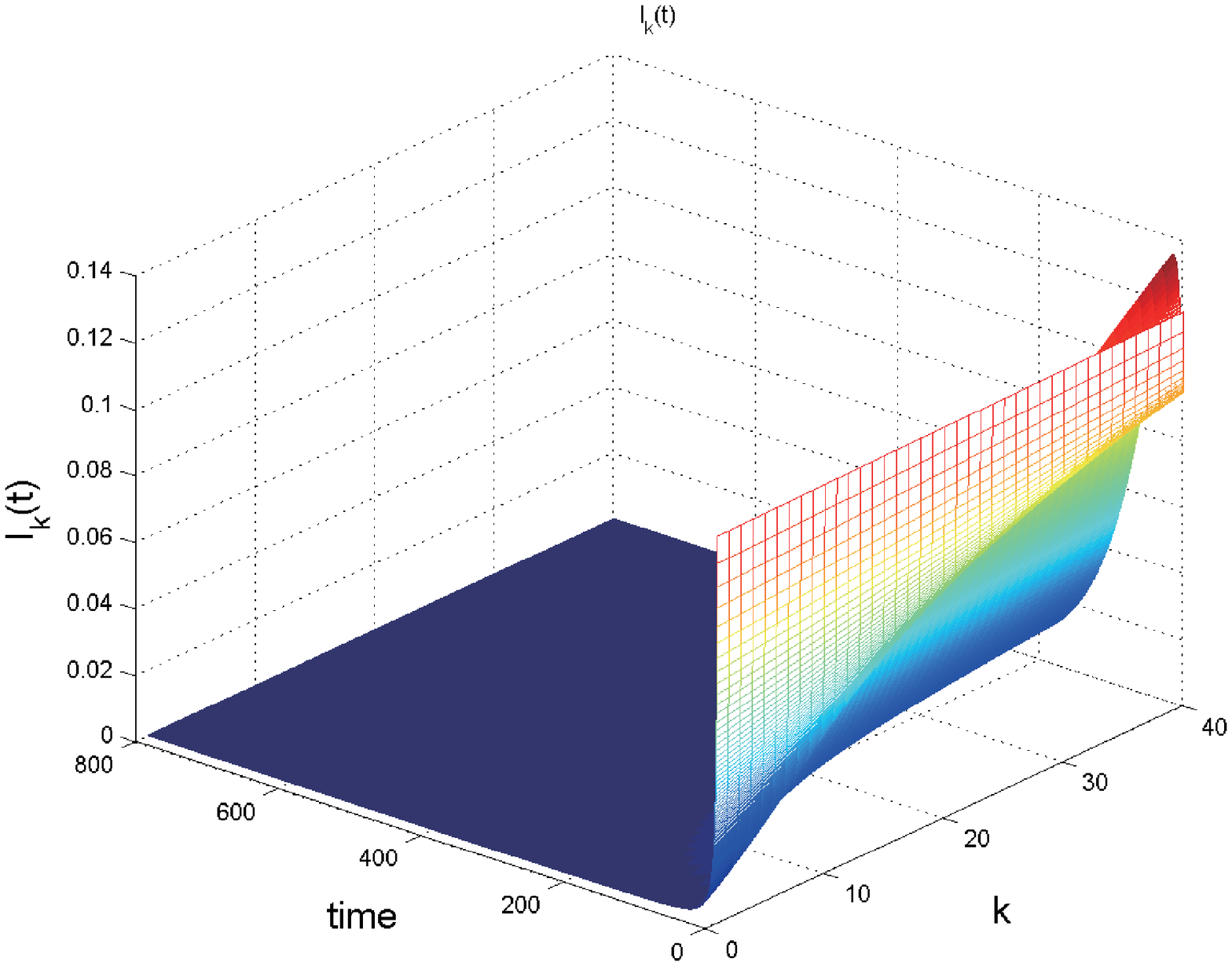}
\label{Fig2(b)}}
\caption{\footnotesize{(Color online) Here $\varphi(k)=k,\beta(\tau)=\frac{\tau(200-\tau)}{15000},\gamma(\tau)=\frac{1}{1+\tau},\mu=0.06,b=0.07,R_{0}=0.6066$}}
\label{Fig2}
\end{figure}

\begin{figure}[htbp]
\centering
\subfigure[Dynamics of I(t) subject to time $t$]{\includegraphics[width=6.3cm]{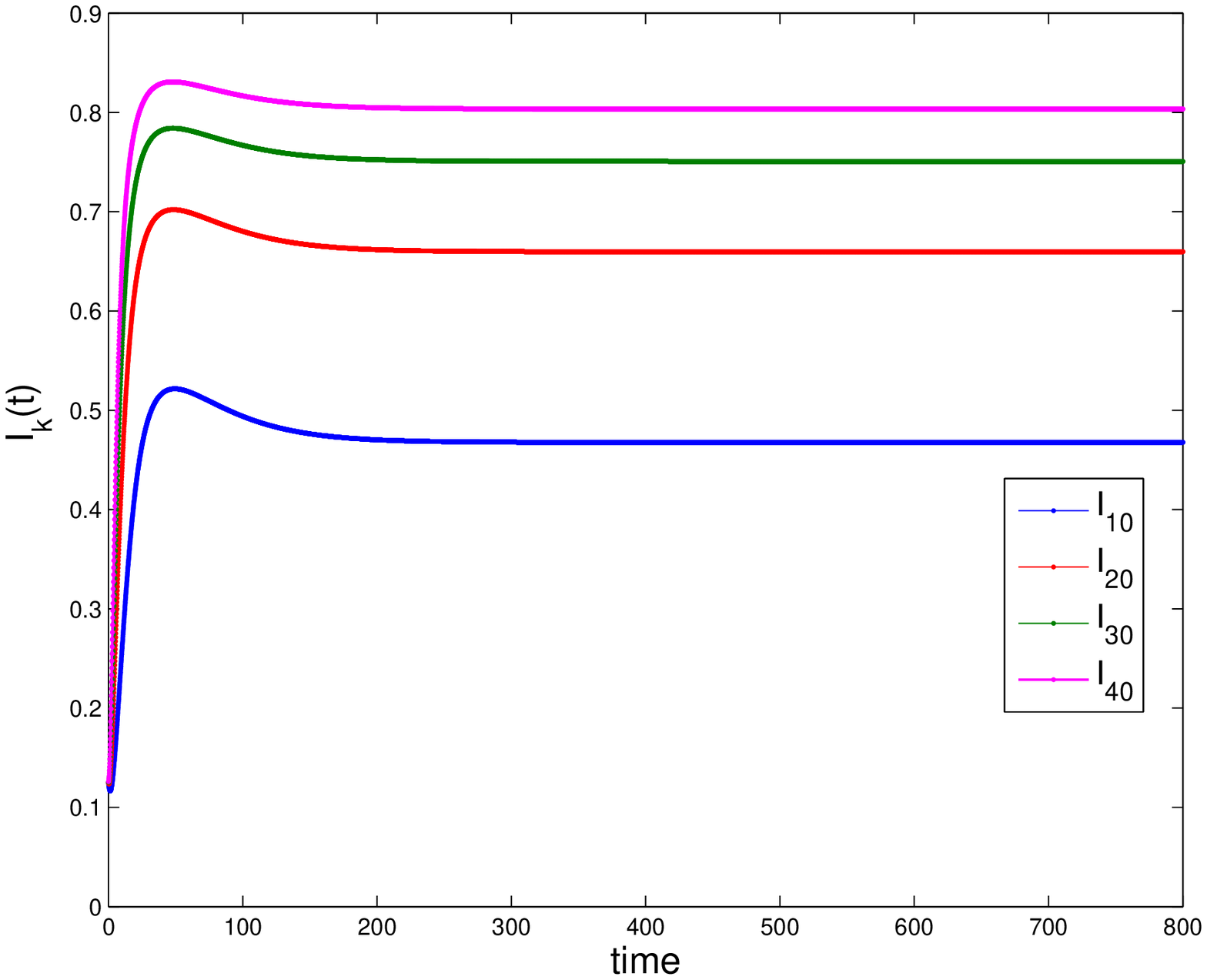}
\label{Fig3(a)}}
\subfigure[Dynamics of I(t) subject to time $t$]{\includegraphics[width=7.2cm]{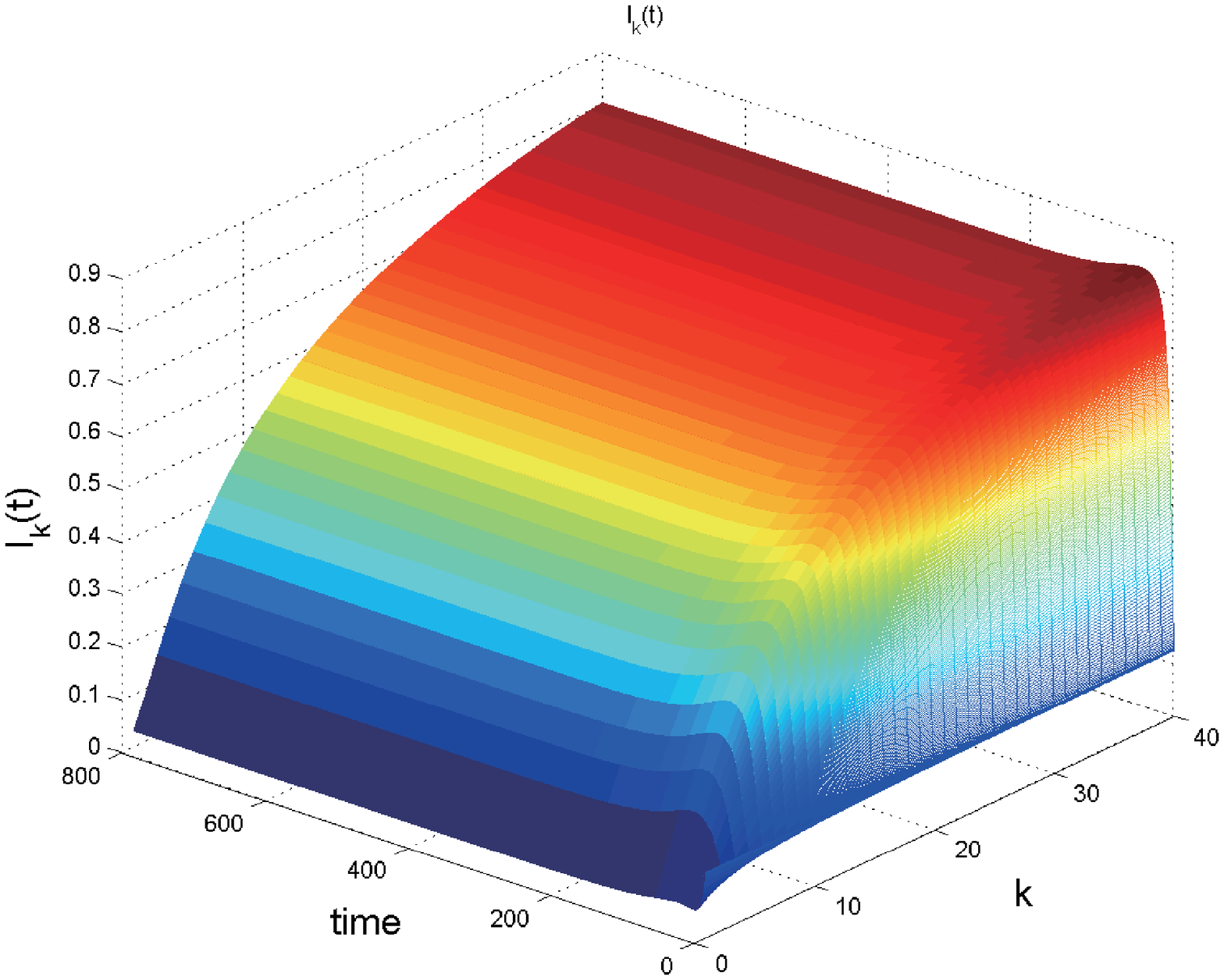}
\label{Fig3(b)}}
\caption{\footnotesize{(Color online) Here $\varphi(k)=k,\beta(\tau)=\frac{\tau(200-\tau)}{15000},\gamma(\tau)=\frac{1}{1+10\tau},\mu=0.06,b=0.07,R_{0}=3.4798$}}
\label{Fig3}
\end{figure}


\section{Conclusions and discussions}

In this paper, we propose and analyze an SIS epidemic model with age-structure on scale-free networks. By using various analytic methods, we demonstrate the asymptotic smoothness of solutions and uniform persistence of the system ~\eqref{eq21} via analyzing its limiting system~\eqref{eq25} with Volterra integral equations. We found that the basic reproduction number is not only related to the network structure, but also depends on parameters which relate to each individual's age.
In addition, we showed that the disease-free equilibrium in model ~\eqref{eq25} is globally asymptotically stable if $R_{0}<1$  by analyzing the corresponding
characteristic equations and applying Fatou's Lemma. Meanwhile, if $R_{0}>1$,  the system has a unique endemic equilibrium, which is globally asymptotically stable by constructing proper Volterra-type Lyapunov functionals. Finally, we performed some simulations under different parameters and $\varphi(k)$, which confirmed our theoretical results.

The factor of age-structure is necessary when study some particular diseases, therefore, in order to analyze the spreading mechanism and dynamical behavior of epidemic diseases more realistically, our findings in this paper are valuable for the further study of age-structured models on complex networks.

\subsection*{Acknowledgment}

This work was jointly supported by the NSFC grants 11331009 and 11572181.
And the paper was partly done while visiting the Center for Mathematical Sciences,
at Huazhong University of Science and Technology, Wuhan, China.



\begin{thebibliography}{99}

 \bibitem{1}
Anderson, R. M., and May, R. M., Infectious Diseases of Humans, Oxford Univ. Press, 1991.

 \bibitem{2}
Bailey, N.T.J., The Mathematical Theory of Infectious Disease, (2nd ed.), New York: Hafner, 1975.

  \bibitem{3}
Kermack, W.O., Mckendrick, A.G., Contributions to the mathematical theory of epidemics.
Proc. Roy. Soc. 1927, A 115: 700-721.

  \bibitem{4}
Kermack, W.O., Mckendrick, A.G., Contributions to the mathematical theory of epidemics.
Proc. Roy. Soc. 1932, A 138: 55-83.

 \bibitem{5}
Levin, S.A., and Durrett, R., From individuals to epidemics.
Philos. Trans. Roy. Soc. London, 1996, B 351: 1615-1621.

 \bibitem{6}
Read, J.M., and Keeling, M.J., Disease evolution on networks: the role of contact structure.
Proc. Roy. Soc., 2003, B 270: 699-708.

 \bibitem{7}
Keeling, M.J., and Eames, K.T.D., Networks and epidemic models.
J. Roy. Soc. Interface, 2005, 2: 295-307.

 \bibitem{8}
Durrett, R., and Jung, P., Two phase transitions for the contact process on small worlds.
Stochastic Proc. \& Their Appl., 2007, 117: 1910-1927.

\bibitem{9}
Pastor-Satorras, R. and Vespignani, A., Epidemic spreading in scale-free networks.
Phys. Rev. Lett., 2001, 86: 3200-3203.

\bibitem{10}
Pastor-Satorras, R., and Vespignani, A., Epidemic dynamics in finite size scale-free networks.
Phys. Rev. E, 2002, 65: 035108.

\bibitem{11}
Pastor-Satorras, R., and Vespignani, A., Epidemic dynamics and endemic states in complex networks.
Phys. Rev. E, 2001, 63: 066117.

\bibitem{12}
Bai, W.J., Zhou, T., et al., Immunization of susceptible-infected model on scale-free networks.
Physica A 2007, 384: 656.

\bibitem{13}
Olinky R., Stone L., Unexpected epidemic threshold in heterogeneous networks: the role of disease transmission. Phys. Rev. E, 2004, 70: 030902.

\bibitem{14}
Zhu, G.H., Fu, X.C., Chen, G.R., Spreading dynamics and global stability of a generalized epidemic model on complex heterogeneous networks.
Appl. Math. Model., 2012, 36: 5808-5817.

\bibitem{15}
Fu, X.C., Small, M., Walker, D.M., Zhang, H.F., Epidemic dynamics on scale-free networks with piecewise linear infectivity and immunization. Phys. Rev. E, 2008, 77: 036113.

\bibitem{16}
Wang, L., and Dai, G., Global stability of virus spreading in complex heterogeneous networks.
SIAM J. Appl. Math., 2008, 68: 1495-1502.

\bibitem{17}
Zhang, J., and Jin, Z., The analysis of an epidemic model on networks.
Appl. Math. Comput., 2011, 217: 7053-7064.

\bibitem{18}
Hoppensteadt, F., An age-dependent epidemic model.
J. Franklin Inst., 1974, 297: 325-338.

\bibitem{19}
Hoppensteadt, F., Mathematical Theories of Populations: Demographics, Genetics and Epiemics, SIAM Publications, Philadelphia, 1975.

\bibitem{20}
Magal, P.,  McCluskey, C.C., Webb, G.F., Lyapunov functional and global asymptotic stability for an infection-age model.  Appl. Anal.  2010, 89: 1109-1140.

\bibitem{21}
Huang, G.,  Liu, X.,  Takeuchi, Y.,
Lyapunov functions and global stability for age-structured HIV infection model.
SIAM J. Appl. Math.,  2012, 72: 25-38.

\bibitem{22}
Tarkhanov, N., Lyapunov stability for an age-structured population model.
Ecol. Model.,  2008, 216: 232-239.

\bibitem{23}
Liu, L.L., Wang, J.L., Liu, X.N., Global stability of an SEIR epidemic model with age-dependent latency and relapse. Nonl. Anal. RWA,  2015, 24: 18-35.

\bibitem{24}
Barab\'{a}si, A.L., Albert, R., Emergence of scaling in random networks.
Science, 1999, 286: 509-512.

\bibitem{25}
Liu, J., Tang, Y., and Yang, Z., The spread of disease with birth and death on networks.
J. Stat. Mech., 2004, 2004: P08008.

\bibitem{26}
Zhu, G.H., Chen, G.R., Xu, X.J., Fu, X.C., Epidemic spreading on contact networks with adaptive weights.
J. Theor. Bio., 2013, 317: 133-139.

\bibitem{27}
Hale, J.K. , Functional Differential Equations, Springer, Berlin, 1971.

\bibitem{28}
Browne, C.J., Plyugin, S.S.,  Global analysis of age-structured within-host virus model.
Discret. Cont. Dyn. Syst. Ser. B,  2013, 18: 1999-2017.

 \bibitem{29}
 Webb, G.F., Theory of Nonlinear Age-Dependent Population Dynamics, Marcel Dekker, New York, Basel, 1985.

\bibitem{30}
Hale, J.K., Asymptotic Behavior of Dissipative Systems, in: Math. Surv. Monogr., vol. 25, Am. Math. Soc., Providence, RI, 1988.

\bibitem{31}
 Adams, R.A.,  Fournier, J.J.F.,  Sobolev Spaces, (2nd ed.),
 Elsevier/Academic Press, Amsterdam, New York, 2003.

\bibitem{32}
 Thieme, H.R., Uniform persistence and permanence for non-autonomous semiflows in population biology.
  Math. Biosci.,  2000, 166: 173-201.

\bibitem{33}
 Hirsch, W.M.,  Hanisch, H.,  Gabriel, J.P.,  Differential equation models of some parasitic infections: methods for the study of asymptotic behavior. Comm. Pure Appl. Math.,  1985, 38: 733-753.

\bibitem{34}
 Li, M.Y., Shuai, Z.S., Global-stability problem for coupled systems of differential equations on networks.
 J. Diff. Eqs., 2010, 248: 1-20.

\bibitem{35}
Hale, J.K., Verduyn, L.S.M., Introduction to Functional Differential Equations, New York:
Appl. Math. Sci., Springer, 1993.

\bibitem{36}
Zhou, T., Liu, J., Bai, W., Chen, G.R., and Wang, B., Behaviors of susceptible-infectedepidemics on scale-free networks with identical infectivity. Phys. Rev. E, 2006, 74: 056109.

\bibitem{37}
Zhang, H.F., and Fu, X.C., Spreading of epidemics on scale-free networks with nonlinearinfectivity.
Nonl. Anal. TMA, 2009, 70: 3273-3278.

\bibitem{38}
Zhu, G.H., Fu, X.C., Chen, G.R., Global attractivity of a network-based epidemic SIS
model with nonlinear infectivity. Commun. Nonl. Sci. Numer. Simul.
2012, 17: 2588-2594.





\end{thebibliography}
\end{document}